\documentclass[12pt]{article}

\usepackage{amsmath}
\usepackage{amsthm}
\usepackage{amssymb}
\usepackage{amstext} 
\usepackage{float}
\usepackage{graphicx}
\usepackage{subfig}
\usepackage{bbm}
\usepackage[pagewise]{lineno}
\newtheorem{theorem}{Theorem}

\newtheorem{remark}[theorem]{Remark}

\usepackage{listings}
\usepackage{nicefrac}
\usepackage{color}
\usepackage{epstopdf}

\usepackage{amssymb}

\usepackage{lineno}

\usepackage{tabularx}
\newcolumntype{+}{>{\global \let \currentrowstyle \relax}}
\newcolumntype{^}{>{\currentrowstyle }}

\usepackage{booktabs}
\usepackage{multirow}

\usepackage{lscape}

\allowdisplaybreaks

\newcommand\BibTeX{{\rmfamily B\kern-.05em \textsc{i\kern-.025em b}\kern-.08em
		T\kern-.1667em\lower.7ex\hbox{E}\kern-.125emX}}

	\title{Disease contagion models coupled to crowd motion and mesh-free simulation}

\author{Parveena Shamim Abdul Salam\\
	\small{Department of Mathematics, TU Kaiserslautern}\\
	\small{Kaiserslautern, 67663,
		Germany} \\
	\small{parveena@mathematik.uni-kl.de }\\[.3cm]
	Wolfgang Bock\\
	\small{Department of Mathematics, TU Kaiserslautern}\\
	\small{Kaiserslautern, 67663,
		Germany} \\
	\small{bock@mathematik.uni-kl.de} \\[.3cm]
	\small{Department of Mathematics, TU Kaiserslautern}\\
	\small{Kaiserslautern, 67663,
		Germany} \\
	\small{	klar@mathematik.uni-kl.de} \\[.3cm]
	\small{Department of Mathematics, TU Kaiserslautern}\\
	\small{Kaiserslautern, 67663,
		Germany} \\
	\small{ tiwari@mathematik.uni-kl.de} \\[.3cm] }	

\begin{document}

\maketitle

	\begin{abstract}
		Modeling and simulation of disease spreading in pedestrian crowds has been recently become a topic of increasing relevance. In this paper, we consider  the influence of the crowd motion in a complex dynamical environment on the course of infection of the pedestrians.
		To model the pedestrian dynamics we consider a kinetic equation for multi-group pedestrian flow based on a  social force model coupled with an Eikonal equation. This model is coupled  with a non-local SEIS contagion model for disease spread, where
		besides the description of local contacts also the influence of contact times has been modelled. Hydrodynamic approximations of the coupled system are derived. 
		Finally, simulations of the hydrodynamic model are carried out using a mesh-free particle method.
		Different numerical test cases   are investigated  including uni- and bi-directional flow in a passage with and without obstacles. 
	\end{abstract}
	
	\textbf{keywords:}{pedestrian flow models; disease spread models; multi-group macroscopic equations;   particle methods}\\
	

\section{Introduction}

%Pedestrian flow modelling  has attracted the interest of a large number of scientists from different research fields, as well as  planners and designers. While planning the architecture of buildings one might be interested in how people move around their intended design so that  shops, entrances, corridors, emergency exits and seating can be placed in useful locations. Pedestrian models are helpful in improving efficiency and safety in public places such as airport terminals, train stations, theatres and shopping malls. They are not only used as a tool for understanding pedestrian dynamics at public places, but also support transportation planners or managers to design timetables. 
The COVID-19 pandemic struck the everyday life worldwide. To avoid further spread in absence of a vaccine or a well-established medical therapy, most countries in the world invoked non-pharmaceutical intervention measures as extensive backtracking, testing and quarantine \cite{ferguson, Chowdhury}. One key point is the contact reduction, which is often based on the minimization of the number of contacts and contact time \cite{MOCOS, ferguson, Walker}. 
In many cases, such as in schools or large working facilities social distancing is not always possible, especially at the end of classes or shifts, crowds are formed to leave the facilities. Simulations of disease spread in a moving crowd could give valuable information about how risky these mass events are and support the design of such intervention measures.

Modelling of crowd motion has been investigated  in many works on different levels of description. Microscopic (individual-based)  models based on Newton type equations as well as vision-based models or cellular automata models and agent-based models have been developed, see Refs. \cite{heMo, helbing2, B01,  deg, maury, piccoli}.  Associated macroscopic  pedestrian flow equations involving equations for density and mean velocity of the flow have been  derived as well and 
investigated thoroughly, see Refs. \cite{bell,  maury, helbing, etikyala2014particle, ama, Col1,Col2}. 
An elegant way to include geometrical information and goal of the pedestrians into these models via the additional solution of an eikonal equation has been developed by Hughes, see Refs.
\cite{ama, di, hughes, hughes1,num}. 
For the derivation and relations of the aboves approaches to each other  we refer to  Refs. \cite{di2,CPT,etikyala2014particle}. 
More complex geometries and obstacles have been included into the models  by many authors, see, for example, Refs.
\cite{piccoli,TGD14,TCP06}.
Pros and Cons of these models have been discussed in various reviews, we refer to  Refs. \cite{bell,HJ10,ABF,BBK13} for a detailed discussion of the different approaches.

On the other hand, there is a vast literature on  disease spread models, see Ref.
\cite{Bellomopand} for an overview on multiscale models and connection to crowd dynamics. Moreover, we refer to 
Refs.  \cite{H,K1,K2,MG} for a small selection of papers on mathematical models based on   dynamical systems and to Refs. \cite{house,perez,Bock} for agent-based and network models. 
Models coupling crowd motion and contagion dynamics are far less investigated. We refer to Ref.  \cite{KQ}
for a recent investigation coupling a crowd motion model with a contagion model in a one-dimensional situation.

One main objective of the present  paper is to include a kinetic disease spread model in the form of a multi-group equation into a kinetic pedestrian dynamics model, derive  hydrodynamic approximations and provide an efficient  numerical simulation of the coupled model for complex two-dimensional geometries.
For the pedestrian flow model  we consider a kinetic equation for multi-group pedestrian flow based on a  social force model coupled with an Eikonal equation  to model geometry and goals of the pedestrians. This model is coupled  with a non-local  contagion model for disease spread, where
local contacts as well  as the influence of contact times is included. 
A second objective is the extension of the methodology to situations with more complex geometries and 
moving objects in the computational domain. This is a way to model, for example, the detailed interaction of pedestrians with larger geometrically extended objects like cars in a  shared space
environment. See \cite{BKM14} for another approach to the  interaction of pedestrians and vehicles.
The numerical simulation is, as in Refs. \cite{etikyala2014particle,KT17} , based on  mesh-free particle methods  \cite{tiwari} for the solution of the Lagrangian form of the hydrodynamic equations. Such a methodology gives an efficient
and elegant way to solve the full coupled problem in a complex environment.

The paper is organized in the following way: in section 2 a kinetic model for pedestrian dynamics with disease spread  is presented. The section contains also  the associated hydrodynamic equations derived from a moment closure approach. The meshfree particle methods used in the simulations is  shortly described in  Section 3. The section contains also the results
of the numerical simulations. The solutions of the macroscopic equations are presented for  different physical situations and parameter values including uni- and bi-directional flow in a two-dimensional passage without obstacles and with fixed and moving obstacles.

\section{Equations}

We use a  kinetic model for the evolution of the distribution functions  of susceptible, exposed and infected pedestrians  as a starting point and derive associated hydrodynamic equations.

\subsection{Kinetic evolution equation}
We consider an equation for the evolution of    pedestrian distribution functions $f^{(k)}= f^{(k)}(x,v,t), k=S,E,I$.
Here $f^S$ stands for the distribution of susceptible pedestrians, $f^E$ for the exposed pedestrians and 
$f^I$ for the infected. 
The evolution equations are  given by
%\begin{align}
%& \partial_t f+v\cdot\nabla_r f+\omega\cdot \nabla_\theta f\notag\\
%&+ \nabla_{v}\left(\left(\gamma(u-v)-\frac{1}{m}\int\int 
%\nabla_{r}U(r,\bar{r},\theta,\bar{\theta})\left(\int\int 
%f(\bar{z})d\bar{v}d\bar{\omega}\right)d\bar{r}d\bar{\theta}\right) \cdot 
%f\right)\notag \\
%&+\nabla_v\left( -\nabla_r V_1(r)\cdot f -\frac{A^2}{2}(vf+\nabla_v f) \right)\notag\\ 
%&+\nabla_{\omega} \left(\left(\bar{\gamma}(g(\theta,u)-\omega)-\frac{1}{I_c}\int\int 
%\nabla_{\theta}U(r,\bar{r},\theta,\bar{\theta})\left(\int\int 
%f(\bar{z})d\bar{v}d\bar{\omega}\right)d\bar{r}d\bar{\theta}\right)\cdot 
%f\right)\notag\\
%&+\nabla_\omega \left( -\nabla_\theta V_2(\theta)\cdot f -\frac{B^2}{2}(\omega f+\nabla_\omega f) \right) =0\ ,
%\label{eq:meanfield}
%\end{align}
%\begin{eqnarray}\label{eq:2.9}
%\frac{dx_i^{(k)}}{dt} &=& v_i^{(k)}, \nonumber \\
%\frac{dv_i^{(k)}}{dt} &=& - \sum_{l=1}^{M}\frac{1}{N_k} \nabla_{x_i^{(k)}} \sum_{i\neq j} U^{(k,l)}(\mid x_i^{(k)} - x_j^{(l)}\mid) + G^{(k)}(x_i,v_i^{(k)},\rho_{i}^{(k)})
%\end{eqnarray} 
\begin{equation}\label{eq:kin}
	\partial_t f^{k} + v \cdot \nabla_x f^{k} + R f^k = T^{k} 
\end{equation}
with $k =S,E,I$.
The operators $R$ and $T^k$ are given by the following definitions.
\begin{equation*}
	\begin{aligned}
		R  f^k (x,v)=\nabla_v \cdot \left([G(x,v;\Phi, \rho) -
		\nabla_x U \star \rho (x)]f^k\right).
	\end{aligned}
\end{equation*}  
with
\[
\rho = \rho^S+\rho^E + \rho^I,
\]
where 
\[
\rho^k(x) = \int f^k(x,v) dv .
\]
Here, $U$ is an interaction potential  describing the local  interaction of the pedestrians and  $\star$ denotes the convolution.  Common choices  for the interaction potential are  purely repelling potentials like spring-damper potentials or attractive-repulsive potentials like the   Morse potential. In this paper we have, for simplicity,  considered a Morse potential without attraction  given by 
\begin{equation}\label{eq:morse}
	U = C_r exp\left( -\frac{|x-y|}{l_r}\right), 
\end{equation}
where $C_r$ is the repulsive strength and $l_r$ is the length scale. 
The part of the forcing term involving $G$ describes the influence of the geometry on the pedestrian's motion and a potential long range interaction between the pedestrians, see below for a detailed description.
The operators $T^k$ are defined using  a SEIS-type kinetic disease spread model leading to 
\begin{eqnarray}\label{eq:SIS}
	T^S = \nu f^I- \beta_I f^S \nonumber  \\
	T^E=  \beta_I f^S - \theta f^E \\
	T^I=  \theta f^E - \nu f^I  \nonumber
\end{eqnarray}
with  constants $\nu, \theta$, see Remark \ref{SEIS} below,  and the non-local  infection rate $\beta_I  = \beta_I(x,v; f^S(\cdot ),f^E(\cdot ),f^I(\cdot )) $
depending in a non-local way on the rate of infected persons, compare Refs. \cite{KQ,Bock} for similiar approaches. We define 
\begin{eqnarray}
	\label{eq:contagion}
	\beta_I=  \int\frac{1}{\rho(y)} \int \phi( x- y, v -w)  f^I (y,w)dw dy .
\end{eqnarray}

The kernel $\phi$ in the infection rate is chosen as
\begin{eqnarray}
	\label{eq:phi}
	\phi =  \phi(x,v )  =  i_o \phi_X(x )  \phi_V(  v  ) .
\end{eqnarray}
with $\int \phi_X (x)d x = 1 = \int \phi_V (v) dv $. Here,  $\phi_X$ is determined as a decaying function of $\vert x \vert$ to take into account the effect that infections between pedestrians are more probable the closer pedestrians  are approaching  each other.  $\phi_V$ is chosen in a similar way depending on $\vert v \vert$ to  take into account the fact, that infections are more probable the longer the interacting  pedestrians stay close to each other, that means the smaller their relative velocities are. The parameter $i_o$ is determined by the infectivity.  We refer to the section on numerical results for the exact definition of these kernels.

Finally, $G$ is given by
\begin{equation}\label{eq:eikonaltec}
	G(x,v;\Phi,\rho)  = \frac{1}{T} \left( - V(\rho(x))\frac{\nabla \Phi(x)}{\Vert\nabla \Phi (x)\Vert} - v\right),
\end{equation}
where $\Phi$ is determined by the coupled solution of  the eikonal equation
\begin{equation}\label{eq:eikonal}V(\rho)\Vert \nabla\Phi\Vert = 1 .
\end{equation}
This describes the tendency of the pedestrians to move with  a velocity given by a speed  $V(\rho)$ and  a direction given by the solution of the eikonal equation. The eikonal equation essentially includes all information about the boundaries
and the desired direction of the pedestrians via 
the boundary conditions. These boundary conditions  for the eikonal equation are chosen in the following way. For walls or for the boundaries of obstacles in the domain we set the value of $\Phi$ at the boundary to a numerically large value.
For ingoing  boundaries, free boundary conditions for $\Phi$ are chosen, whereas for  outgoing boundaries, where the 
pedestrians aim to go, we set $\Phi =0$.

We note that on the one hand, the eikonal equation 
includes the  geometrical information via boundary conditions. On the other hand,  it models   a global reaction of the pedestrians  to avoid  regions of dense crowds via the term $V(\rho$) in (\ref{eq:eikonal}).

\begin{remark}
	\label{proxemics}
	The parameters in the above formulas, in particular in the definition of  (\ref{eq:eikonal}) and (\ref{eq:morse}) have to be chosen  consistent
	with  empirical data, see \cite{BGO,KGT09}.
\end{remark}

%
%
%\begin{remark}
%	In the definition of the acceleration towards the desired  direction (\ref{eq:G}), the speed with which the pedestrians are moving depends on the density around a pedestrian. In certain situations this could lead to unphysical effects, for example, if the pedestrian is approached from behind. A determination of the density including a "vision cone"  could be used here at the expense of a more complicated model.
%\end{remark}
%\begin{remark}
%	A further remark on the above microscopic model concerns the role of the interactions between the pedestrians.
%	Interactions are not only modelled by the interaction potential  $U$, but also by the Hughes type term (\ref{eq:G}).The motivation for the present way of modelling is a distinction between a  short scale 
%	interaction between the pedestrians in direct encounter described by the interaction potential $U$
%	and a  reaction of the pedestrian on a much larger spatial scale on the global density $\rho$  
%	via the solution of the eikonal equation as in the Hughes approach.
%	In the present model, as in the Hughes model,  a knowledge of the density in the whole domain is assumed
%	for this second kind of interaction.
%	This could be changed to certain subregions of the computational domain by restricting the solution of the eikonal equation to these regions.
%\end{remark}
\begin{remark}
	Instead of the social force model used here, one could as well use more sophisticated interaction models,
	see for example \cite{deg,Bailo,Mahato}. We note that the differences between these  models   in the present hydrodynamic context are  small. The behaviour  of the solutions is rather dependent on the choice of the parameters.
\end{remark}

\begin{remark}
	\label{SEIS}
	The dynamical system called the SEIS-model with constant infection rate is given by 
	\begin{align*}
		\frac{dS}{dt }= - \beta IS  + \nu I\\
		\frac{dE}{dt }=  \beta IS -  \theta E \\
		\frac{dI}{dt }=  \theta E -  \nu I
	\end{align*}
	where $\beta$ is the infection rate, $\nu$ the recovery rate and $\theta$  the rate with which exposed persons are becoming infected. Pedestrians are potentially becoming exposed, when they are in contact with infected pedestrians.
	However, exposed pedestrians are only becoming infectious with a certain rate $\theta$. Exposed pedestrians  do not infect
	other pedestrians. Usually, in the situations  and on the time scales under consideration here, $\nu$  and $\theta$ are very small 
	and set to zero in  numerical simulations such that the number of infected pedestrians remains constant during the simulation.
\end{remark}

\begin{remark}
	A key difference between other agent-based models as e.g.~\cite{Bock} is the dependence of the infection on the relative velocity of the agents.  This has  direct implications on the process of infection during the dynamics, as can be seen the the case studies.
\end{remark}
\subsection{The multi-group hydrodynamic model}
Integrating the kinetic equation against $dv$ and $vdv$ and using a mono-kinetic distribution function to close the resulting balance equations, i.e. approximate 
$$ f^k \sim \rho^k (x)\delta_{u(x)}(v)$$
one obtains a continuity equation for group $k$
\begin{equation}\label{eq:cont}
	\partial_t \rho^k + \nabla_x \cdot (\rho^k u) = \int T^k d v 
\end{equation}
and
the momentum equation
\begin{equation}\label{eq:mom}
	\partial_t u + (u \cdot \nabla_x)u = G(x,u;\Phi,\rho) -  \nabla_x U\star\rho
\end{equation}
with 
the total density $\rho$ given as 
\[
\rho = \rho^S+\rho^E + \rho^I.
\]
Moreover,  
\begin{equation*}
	G(x,u;\Phi,\rho) = \frac{1}{T} \left( - V(\rho(x))\frac{\nabla \Phi(x)}{\Vert\nabla \Phi(x)\Vert} - u\right),
\end{equation*}
where $\Phi$ is determined by solving the eikonal equation
$$ V(\rho)\Vert \nabla\Phi\Vert = 1 .$$

The continuity equations are explicitly written as 
\begin{eqnarray}\label{cont3}
	\partial_t \rho^{S} +  \nabla_x \cdot (u \rho^{S} ) = \nu \rho^I- \beta_I \rho^S  \nonumber\\
	\partial_t \rho^{E} +  \nabla_x \cdot (u \rho^{E} )=  \beta_I   \rho^S - \theta \rho^E\\
	\partial_t \rho^{I} +  \nabla_x \cdot (u \rho^{I} )=  \theta \rho^E - \nu \rho^I. \nonumber
\end{eqnarray}
with $\beta_I = \beta_I(x; \rho^S(\cdot),\rho^E(\cdot), \rho^I(\cdot) , u (\cdot))$, where now
\begin{eqnarray}\label{eq:hydinfrate}
	\beta_I =  \int \phi( x-y ,  u(x)-u(y))   \frac{\rho_I(y)}{\rho(y)} dy .
\end{eqnarray}
%and
%Then compute only 
%\begin{eqnarray}\label{eq:2.5}
%\partial_t \rho + \nabla_x \cdot (\rho u )= 0
%\end{eqnarray}
%\begin{remark}
%	\label{localization}
%	Localisation  
%	\begin{eqnarray}\label{eq:2.5}
%	\beta_I =  \int \phi(x-y ,  \nabla_x u(x))  \alpha_I(y)dy 
%	\end{eqnarray}
%\end{remark}
\begin{remark}
	Multi-group models for pedestrian flows have been also used in different contexts. For example in \cite{Mahato2} a  hydrodynamic multi-group model for pedestrian dynamics with groups of different sizes has been developed and analysed in \cite{Mahato2}.
\end{remark}
\subsection{The hydrodynamic model using volume fractions}

For numerical computations this is rewritten using volume fractions, that means we  solve the continuity equation 
\begin{eqnarray}\label{eq:conttot}
	\partial_t \rho + \nabla_x \cdot (\rho u )= 0
\end{eqnarray}
and the momentum equation 
\begin{equation}\label{eq:momtot}
	\partial_t u + (u\cdot \nabla_x)u = G(x,u;\Phi,\rho) -  \nabla_x U\star\rho.
\end{equation}

Then we compute the volume fractions
$
\alpha^S,\alpha^E,\alpha^I
$
with
$
\alpha^S+\alpha^E+\alpha^I=1
$
as
\begin{eqnarray}\label{eq:fraction}
	\partial_t \alpha^{S} + u \cdot \nabla_x\alpha^{S}  = \nu \alpha^I- \beta^I \alpha^S \nonumber \\
	\partial_t \alpha^{E} + u \cdot \nabla_x  \alpha^{E} =  \beta^I\alpha^S - \theta \alpha^E\\
	\partial_t \alpha^{I} + u \cdot \nabla_x  \alpha^{I} = \theta \alpha^E - \nu \alpha^I \nonumber
\end{eqnarray}
with $\beta_I = \beta_I(x; \alpha^I(\cdot) , u (\cdot))$ defined by
\begin{eqnarray}\label{eq:inffrac}
	\beta_I =  \int \phi( x-y ,  u(x)-u(y))  \alpha_I(y)dy .
\end{eqnarray}
Finally, one computes
$$
\rho^I = \alpha^I\rho, \rho^S = \alpha^S \rho, \rho^E = \alpha^E \rho.
$$

\subsection{Dynamic geometries}

We allow the domain on which the above equations are defined to depend on time. 
In particular, we consider moving  obstacles, which  change their path and speed in order to avoid collisions with the  pedestrians, while moving towards a  specific target. The interaction between the pedestrians and the moving obstacle is additionally modeled by kinematic equations using a repulsive potential similar to the pedestrian-pedestrian interactions. A second Eikonal equation is integrated for modelling the  path of the moving obstacle in the geometry  and the desired destination of the obstacle. The velocity update equation for the obstacle is given as
\begin{eqnarray}
	\frac{dx^O}{dt}& = & v^O  \, , \nonumber \\
	\frac{dv^O}{dt} & = &  -\nabla_x U^O*\rho (x^O )+ G^O(x^O,v^O;\Phi^O, \rho)   \, ,
	\label{Model_obstmicro} 
\end{eqnarray}
where $x^O$ and $v^O$ are the position and velocity of the   obstacle's centre of mass,   $\rho$ is the  density of pedestrians at $x^O$. $U^O$ is an interaction  potential  describing the interaction of the pedestrians  on the obstacle.

$G^O$ is obtained from the gradient of the Eikonal solution $\phi_O$ for the obstacle as  
%with,
\begin{equation*}
	G^O(x^O,v^O;\Phi^O,\rho) = - \frac{1}{T^O} \left( -V^O(\rho(x^O))\frac{\nabla \phi^O(x^O)}{||\nabla \phi^O (x^O)||} - v^O \right) \, , \\
	||\nabla \phi^O|| =  1 .
	%\label{Eikonal_obst}
\end{equation*}

That means, the eikonal equation is for the obstacle only used to include the geometrical informations and the goal of the obstacle.
We note that the action of the obstacle on the pedestrians is given by the solution of equation (\ref{eq:eikonal})
via the boundary conditions at the obstacle's boundaries. In contrast, the action of the pedestrians on the obstacle is given via $U^O$.

%
%\subsection{Microscopic model}
%	An associated microscopic system is given in the following way. First determine the local density by counting  
%	neigbouring particles 
%	$$
%	\rho_i = \sum_j  \delta_S  ( x_i-x_j  )
%	$$
%	with $\int \delta_S (x) dx = 1$. Then, the equations are 
%	\begin{eqnarray}
%	\frac{d x_i}{d t} \,& = & \, v_i \,  \\
%	\frac{d v_i}{d t} \,& = & \, G(x_i,v_i; \Phi,\rho)\,- \sum_{j } \nabla U (x_i - x_j)  \nonumber 
%	\end{eqnarray}
%	together with 
%	\begin{eqnarray}\label{eq:microfrac}
%	\frac{d  \alpha^{S} _i}{d t}   = \nu \alpha^I_i- \beta^I_i \alpha^S_i  \nonumber \\
%	\frac{d \alpha^{E}_i}{dt}  =  \beta^I_i\alpha^S_i - \theta \alpha^E_i\\
%	\frac{d \alpha^{I}_i}{dt} = \theta \alpha^E_i - \nu \alpha^I_i \nonumber
%	\end{eqnarray}
%	with
%	$$ \beta_i^I = \sum_{j \in I}\frac{1}{\rho_j} \phi( x_i - x_j, v_i - v_j ) , $$
%	where the summation is over all infected pedestrians. 
%	\begin{eqnarray}\label{num3}
%	\frac{d x^O}{d t} \,& = & \, v^O \, , \nonumber \\
%	\frac{d v^O}{d t} \,& = & \, G^O(x^O,v^O;\Phi^O,\rho) \,- \, \sum_{j } \nabla U^O (x^O - x_j).
%	\label{Model_Lagrangian}
%	\end{eqnarray}	
%
%

\section{Numerical method and results}

For the  numerical simulation we use a meshfree particle method, which is based on least square approximations.~\cite{fpm,KT17} A Lagrangian formulation of the hydrodynamic equations is used and  coupled to the SEIS model and  the obstacle's kinematic equations Eq.~(\ref{Model_obstmicro}). 

%  We need to exchange the information between these two grid points. For example, we need the solution or the gradient of the eikonal solution to solve the hydrodynamic equations and the density of the hydrodynamic equations is required to solve the eikonal equation. The information from one grid points to another grid points are obtained from the least squares approximation. 

\subsection{Lagrangian equations}
The spatially discretizted  system in Lagrangian form is given by
\begin{eqnarray}
	\label{num1}
	\frac{d x_i}{d t} \,& = & \, u_i \, , \nonumber \\
	\frac{d \rho_i}{d t} \,& = & \, - \rho_i \, \nabla_x  \cdot u_i \, ,\nonumber \\
	\frac{d u_i}{d t} \,& = & \, G(x_i,u_i; \Phi,\rho)\,- \sum_{j } \nabla U (x_i - x_j)  \rho_j dV_j \, , \nonumber \\
\end{eqnarray}
and
\begin{eqnarray}\label{num2}
	\frac{d  \alpha^{S} _i}{d t}   = \nu \alpha^I_i- \beta^I_i \alpha^S_i  \\
	\frac{d \alpha^{E}_i}{dt}  =  \beta^I_i\alpha^S_i - \theta \alpha^E_i\\
	\frac{d \alpha^{I}_i}{dt} = \theta \alpha^E_i - \nu \alpha^I_i
\end{eqnarray}
with
$$ \beta_i^I = \sum_{j } \phi( x_i - x_j,  u_i - u_j ) \, \alpha^I_j \, dV_j ,$$
and
\begin{eqnarray}\label{num3}
	\frac{d x^O}{d t} \,& = & \, v^O \, , \nonumber \\
	\frac{d v^O}{d t} \,& = & \, G^O(x^O,v^O;\Phi^O,\rho) \,- \, \sum_{j } \nabla U^O (x^O - x_j) \rho_j dV_j .
	\label{Model_Lagrangian}
\end{eqnarray}

Here $dV_j$ is the local area around a particle.

\begin{remark}
	Although the equations in Lagrangian form look similiar to a microscopic problem, there are important differences.
	In particular, the 	efficient solution of the continuity equation (compared to a determination of the density in a purely microscopic simulation) and  the consequent availablity of $\rho_i$ on each grid-point allows the efficient use  of the  density, used 
	at various places in the model.
	
\end{remark}

\subsection{Numerical method}

For a description of the mesh-free method for the  pedestrian flow equations in a  fixed rectangular domain, we refer to \cite{tiwari,KT17}.
There, the eikonal equation has been solved on a separate regular mesh on the entire domain using the fast-marching method.  The information from the irregular point cloud, used for solving the fluid equations has been  interpolated to the  eikonal grid
and vice versa. Such an approach requires for a moving obstacle to take special care of the grid points being overlapped by the obstacles. Using an immersed boundary method, they are activated and deacitvated
depending on the location of the obstacle. 

Here, we have followed a different approach. The eikonal equation is directly solved on the 
irregular point cloud used for the hydrodynamic equations.
Thus,  in each time step an eikonal equation on an unstructured grid has to be solved, see \cite{sethian} for the Fast-marching method in this case.
The approximation of the spatial derivatives in the eikonal equation is obtained as for the hydrodynamic equations:   the spatial derivatives at an arbitrary grid point are computed from the values at  its surrounding neighboring grid points using a weighted least squares method. We refer to \cite{eikonal} for a numerical study of the accuracy and complexity of such a method for the eikonal equation. 

%If there is an obstacle in the computational domain, it will overlap the eikonal grids. Those overlapping grids are called as non-active grids. The non-active grids are removed from the neighbor list while solving the eikonal equations. When the obstacle moves we have to activate and diactivate grid points in every time step.
%for regular eikonal grids?
%

Finally, we note that for  a uni-directional flow of pedestrians, we have to solve only one eikonal equation. If there is bi-directional flow, we have to solve an eikonal equation for each direction. The same would be true for several obstacles with different goals.

\subsection{Numerical results}
\begin{figure}
	\centering
	\includegraphics[keepaspectratio=true, angle=90, width=0.45\textwidth]{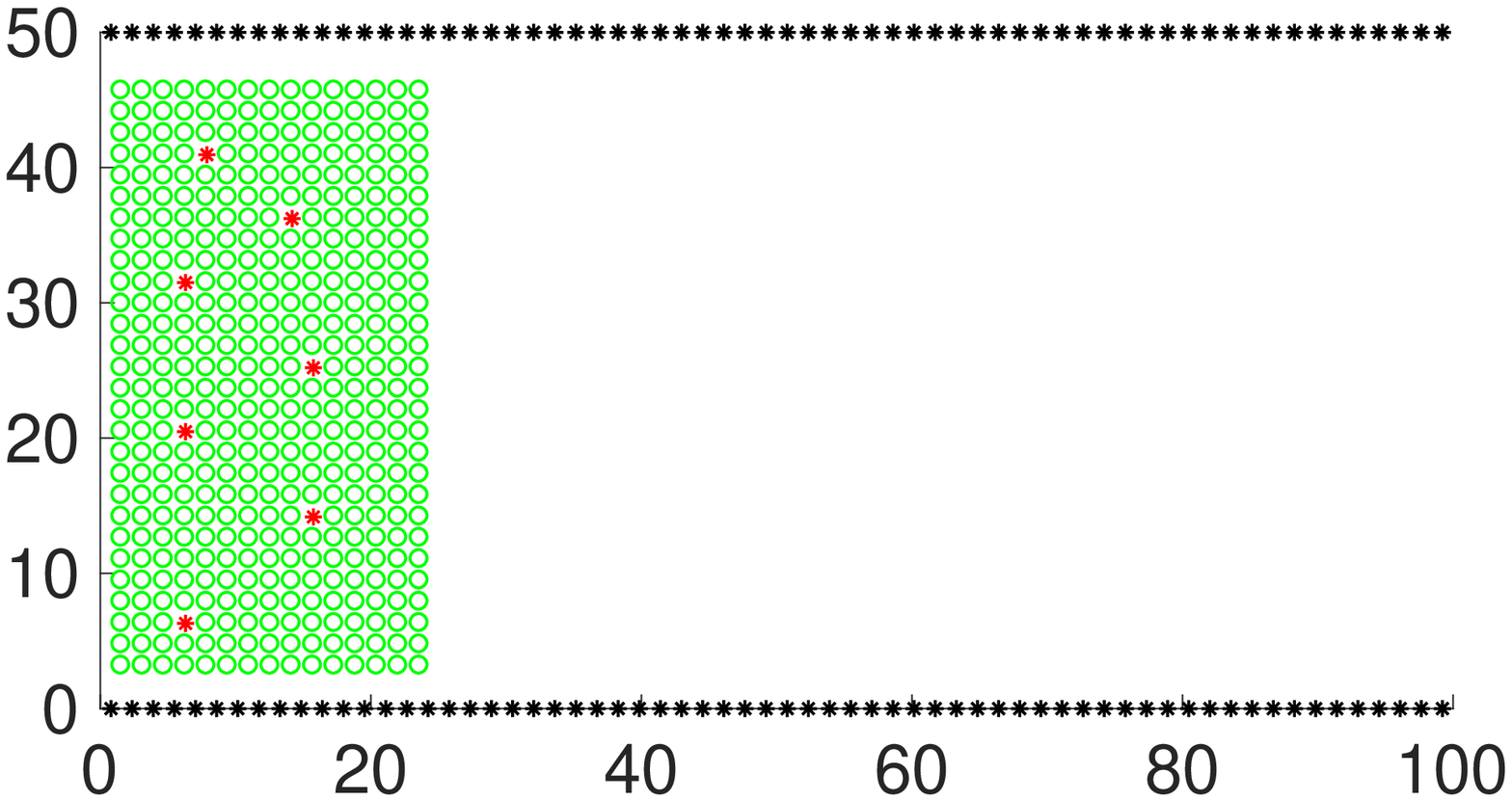}
	\includegraphics[keepaspectratio=true, angle=90, width=0.45\textwidth]{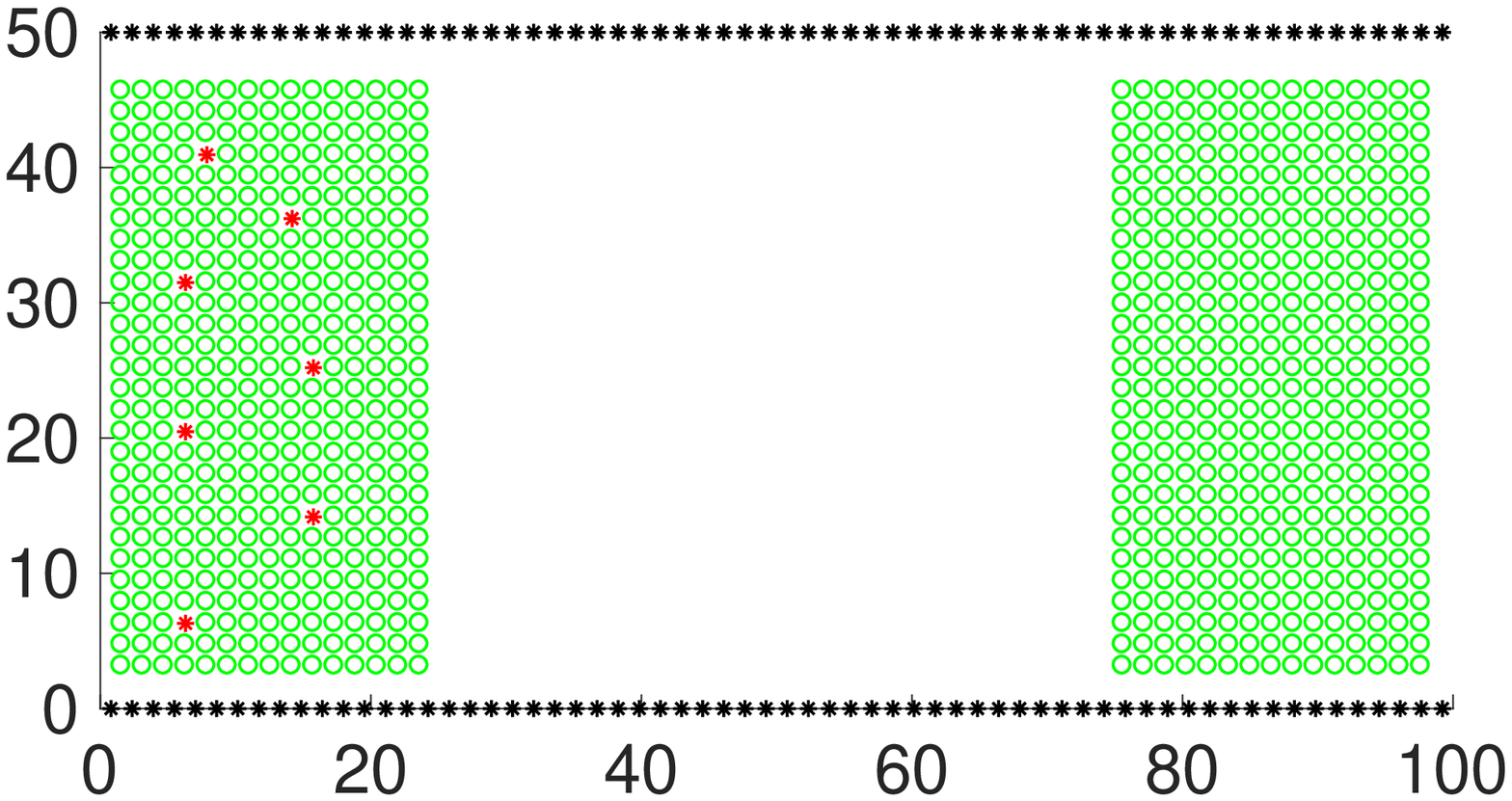}	
	\\ %\vspace{-3cm}
	\includegraphics[keepaspectratio=true, angle=90, width=0.45\textwidth]{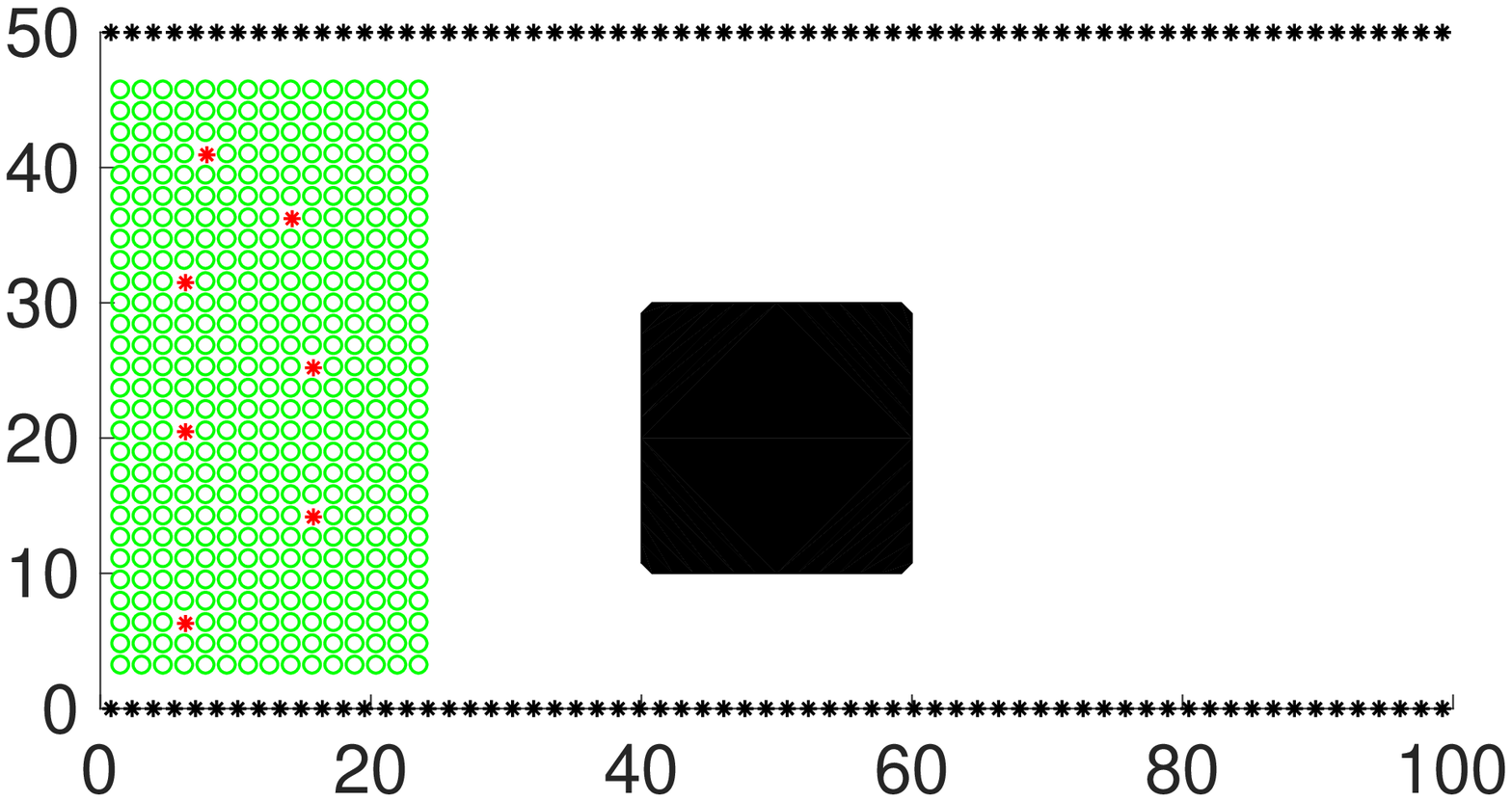}
	\includegraphics[=true, angle=90, width=0.45\textwidth]{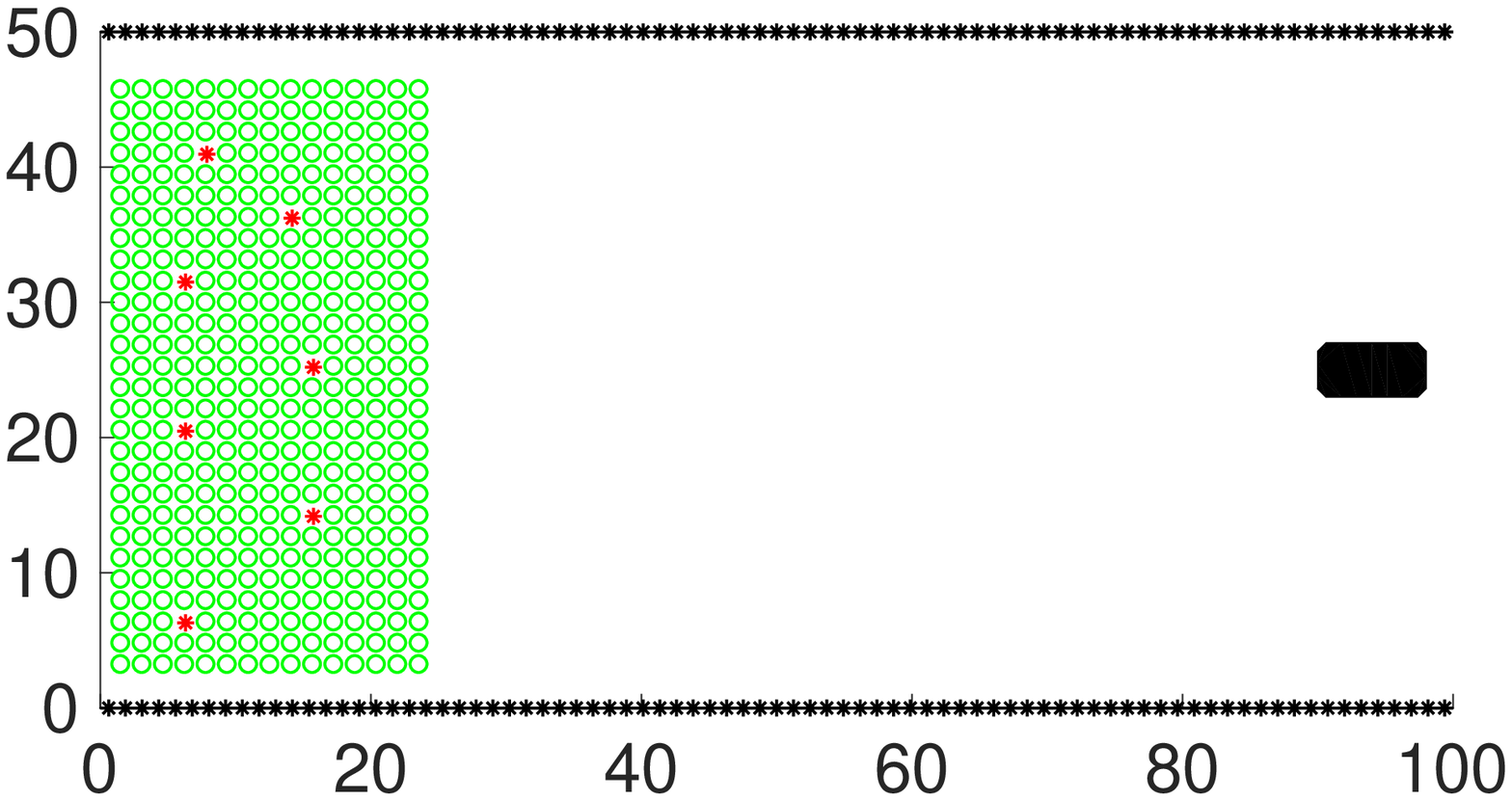}	
	
	%		\includegraphics[keepaspectratio=true, angle=0, width=0.45\textwidth]{figures/fig_with_obstacle_bi_direction_phi_V_1/with_obstacle_bi_direction_Phi_V_1_t_0}	
	%	\vspace{-2cm}	 
	\caption{Initial situation at  $t = 0$. Top row: Uni-directional (left) and bi-directional (right) flow. Bottom row: Fixed obstacle (left) and moving obstacle (right). Red indicate infected, green indicate susceptible pedestrians.}
	\label{phi_v_ne_1_t0}
	\centering
\end{figure}

We have performed numerical simulations of equations (Eq.~(\ref{eq:conttot}) to Eq.~(\ref{eq:inffrac}) and
(\ref{Model_obstmicro})  for different scenarios. In all our simulations, we consider a computational  domain given by a platform or corridor of size $100~m\times 50~m$.
The top and bottom boundaries are rigid walls without any entry or exit. Right and left boundary are  exits for pedestrians and obstacles depending on the situation under consideration. We consider uni- as well as bi-directional flow of the pedestrians. Initially, the pedestrians are distributed as shown in Figure \ref{phi_v_ne_1_t0} with a distance of $1.575 m$ from each other. We consider the cases with and without  obstacles, which are 
either fixed or moving.
We initialize the  infected pedestrians (colored in red) with $\alpha_I=1, \alpha_S=0, \alpha_E=0$, the susceptibles 
(in green) by $\alpha_I=0, \alpha_S=1, \alpha_E=0$.

%  \begin{figure}
% 	\centering
% 	\captionsetup[subfigure]{margin=5pt} %parskip=0pt, hangindent=0pt, indention=0pt, singlelinecheck=true}
% 	\subfloat[microscopic]{
% 		\includegraphics[keepaspectratio=true, width=.5\textwidth]{fig_2d_ped/naive_h_0dot4_R_0dot4.pdf}
% 	} \subfloat[ multi-scale]{
% 		\includegraphics[keepaspectratio=true, width=.5\textwidth]{fig_2d_ped/multiscale_h_0dot4_R_0dot4.pdf}
% 	} 
% 	\caption{ Density plot determined from local limit equation (\ref{pedlocalhydro})  and nonlocal equations (\ref{pednonlocalhydro}) with microscopic and multi-scale approximation 
% 		at $t=12,5$ for $\Delta x = 0.2$ and $R=0.4$.}
% 	\label{compare_ped2}
% \end{figure}
% 

We have used for the infection rate $\beta_I$ the functions 
$\phi_X = \exp(-|x-y|^4)$ and  $\phi_V = \exp(-|u-v|^6)$. 
Moreover,  we choose 
$ V(\rho) = V_{max} \left( 1 - \frac{\rho}{\rho_{max}}\right)$
and 
$ V^O(\rho) $ in the same way with $V_{max}$ substitued by $ V_{max}^O $.
For the parameters  we have used the values given in Table \ref{tbl:comp}.

\begin{table}[ht]
	\center
	{\begin{tabular}{@{}c|c||c|c@{}} \toprule
			Variable & Value & Variable & Value \\ \midrule
			$V_{max}$ & 2 $m/s$ & $\rho_{max}$ & 10 ped/$m^2$ \\
			$V_{max}^O$ & 3 $m/s$ & T & 0.001 $s$ \\
			$C_r=C_r^O$ & 50 & $l_r$ & 2 m\\
			$l_r^{O}$ & 1 m& $i_o$ & 0.04 $m^2/s^2$\\ \bottomrule 
		\end{tabular}}
		\caption{Numerical parameters.}
		\label{tbl:comp}
	\end{table}

	%If the infected pedestrain is approaching to the suspected one, we have used the infection ratio $2 i_0$, which can be computed if $(x-y)\cdot(u-v) < 0$, where $x,y$ are position vectiors and $u,v$ are their velocities. 
	
	During the evolution, 
	infected pedestrians are  colored in red, susceptibles in green and  exposed pedestrians in blue
	according to the values of $\alpha^k$.
	If $\alpha^E > 0.05$ the  colour is switching from green to blue, meaning that the probability of being exposed has exceeded a certain threshold. The red pedestrians remain red throughout the simulations,
	since the recovery rate $ \nu$ is set equal to $0$ in the simulations.  Moreover, since $ \theta $ is also set to $0$  exposed patients are not becoming infected and cannot infect others  in the simulations.

	The fixed  and the moving obstacle considered are  rectangular in shape and initially located as shown in Figure \ref{phi_v_ne_1_t0}.
	
	Explicit time integration  of the equations in Lagrangian form is done with a fixed time step size of $0.001$ in our simulations.

	\subsection {Test-case 1:  fixed obstacle}     
	
	In this first test case we have considered a fixed obstacle and compared  the results to a situation without obstacle.
	We consider bi-directional flow without an obstacle  and uni- and bi-directional flow with an obstacle.
	This is done for the case $\phi_v=0$ (no influence of contact time) and the case where $\phi_V$ is chosen as defined above, that means for a situation where the influence of the contact time is included. 
	The present initial  configuration is chosen in such a way, that  there is no increase of the  number of probably exposed pedestrians,
	if a uni-directional flow without obstacle is considered with or without influence of the contact time.
	
\begin{figure}
		\centering
		\includegraphics[keepaspectratio=true, angle=90, width=0.45\textwidth]{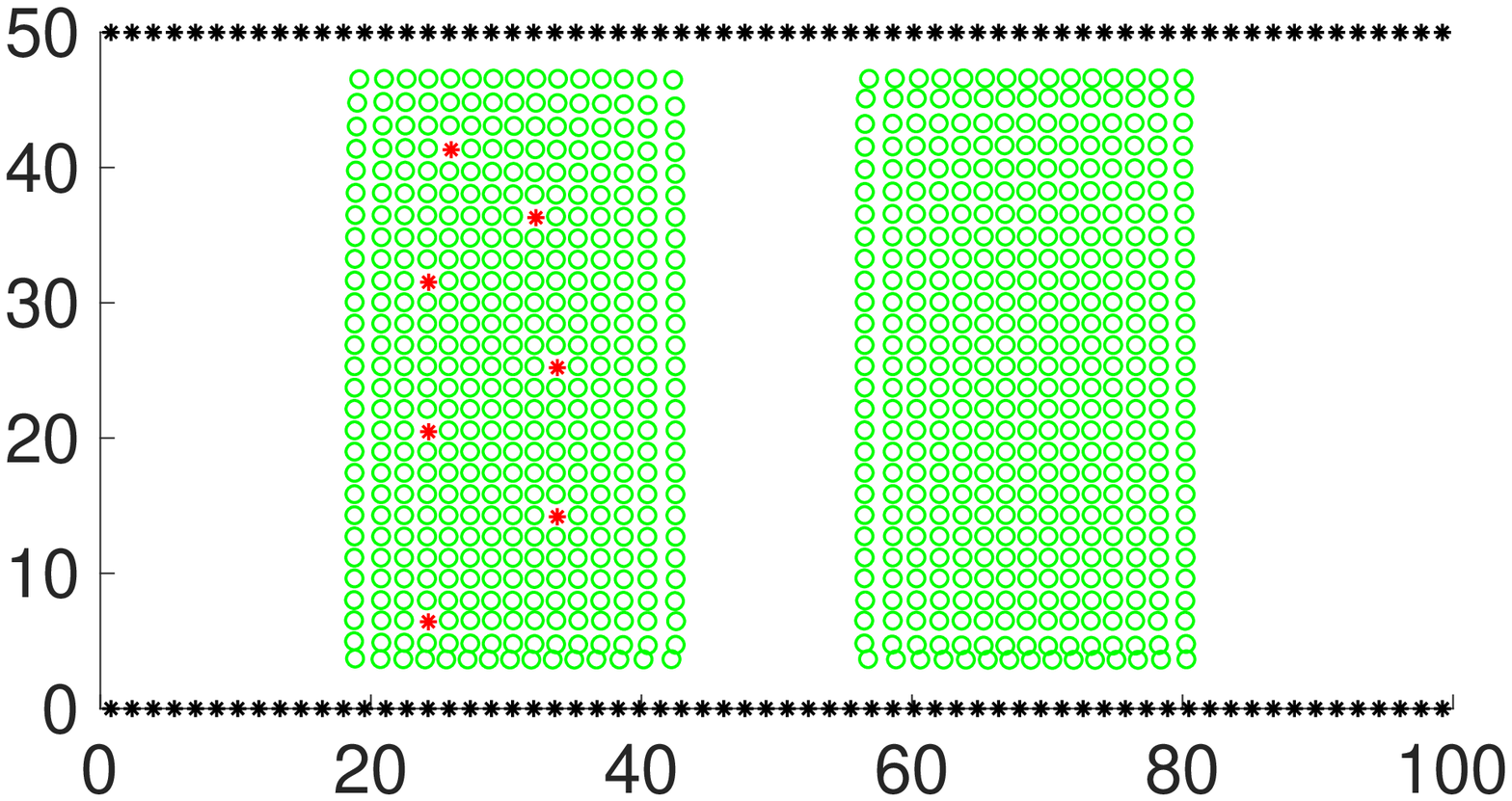}
		\includegraphics[keepaspectratio=true, angle=90, width=0.45\textwidth]{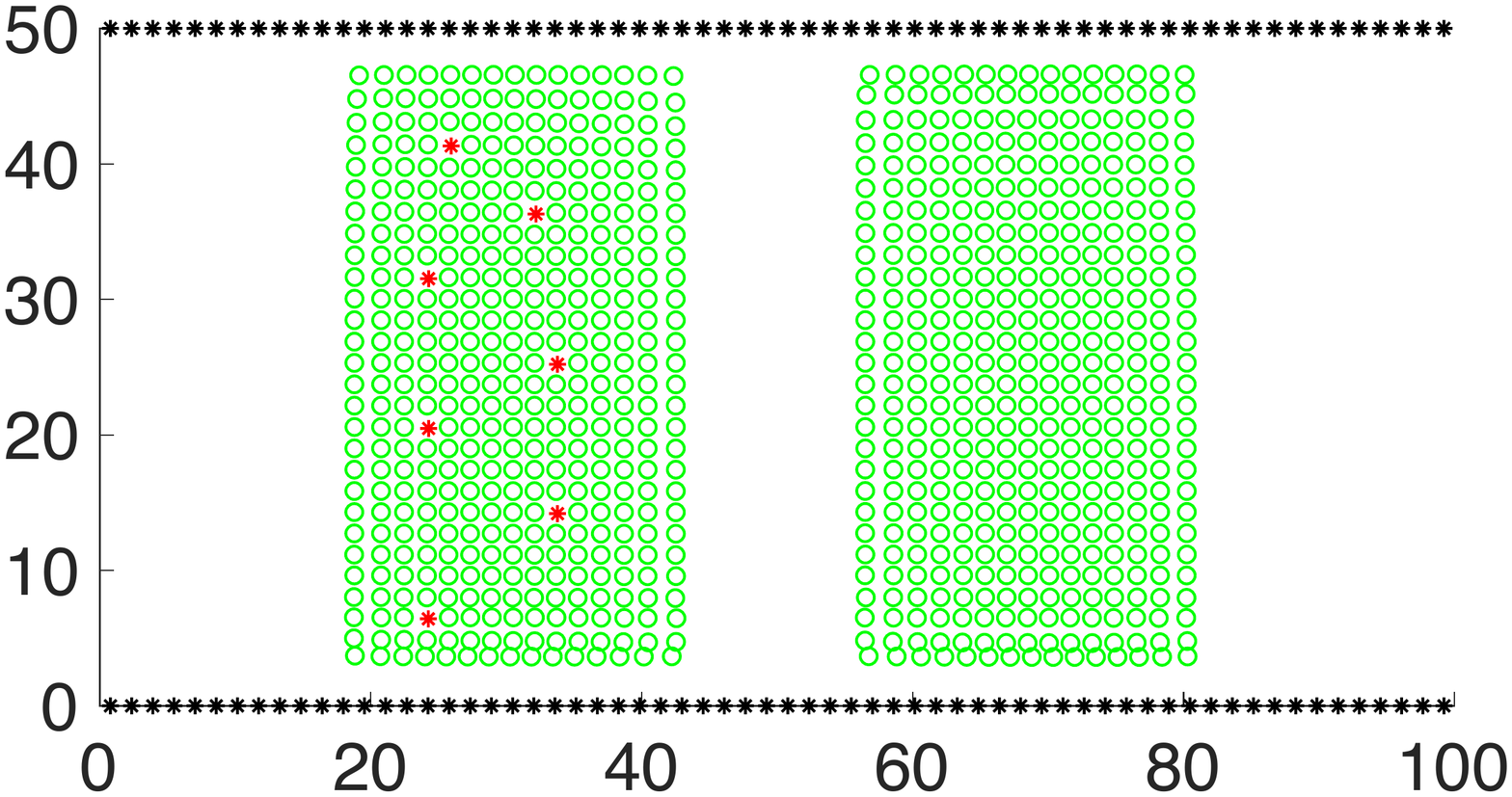}\\
		%	\vspace{-4cm}	
		\includegraphics[keepaspectratio=true, angle=90, width=0.45\textwidth]{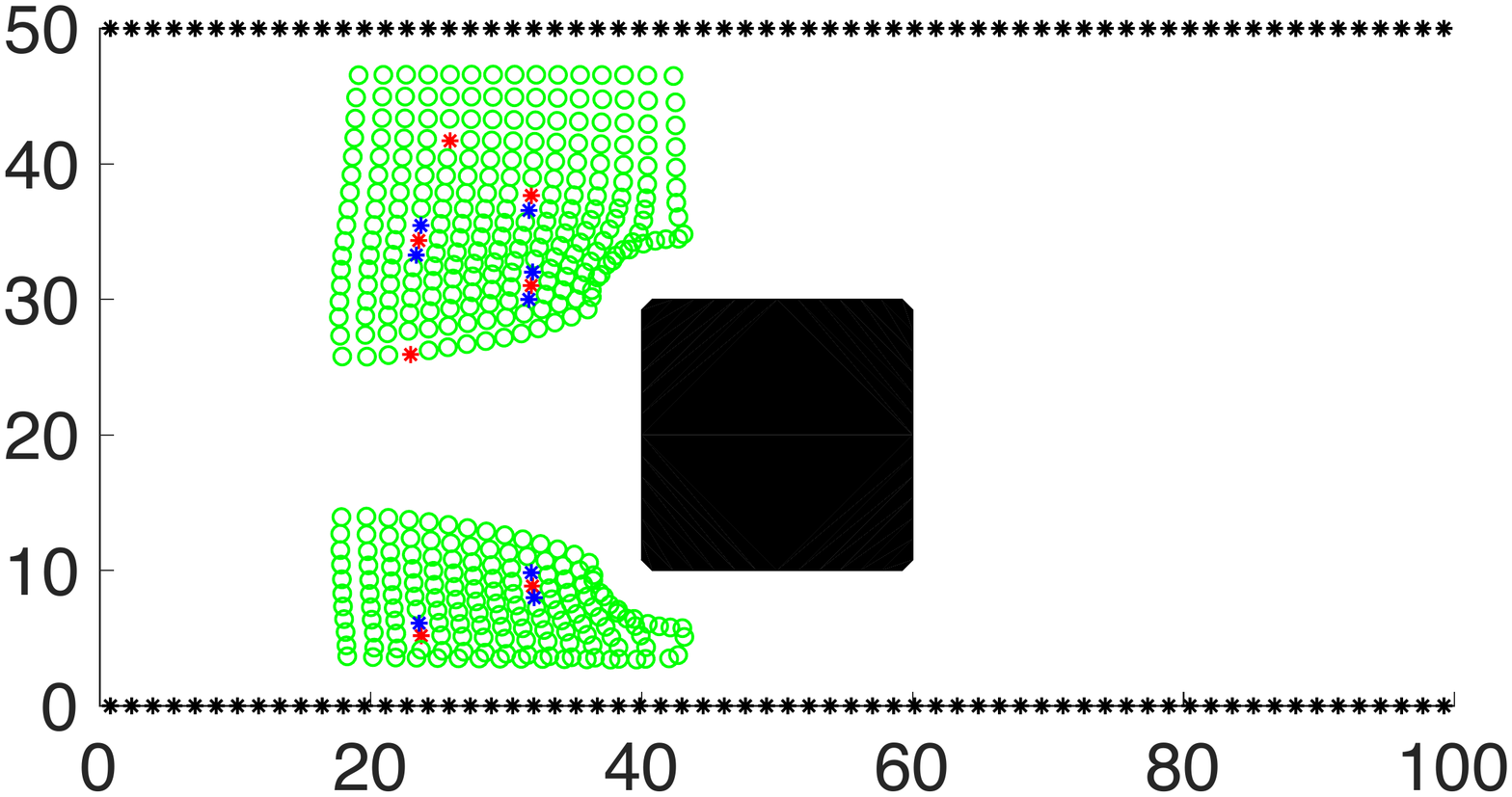}
		\includegraphics[keepaspectratio=true, angle=90, width=0.45\textwidth]{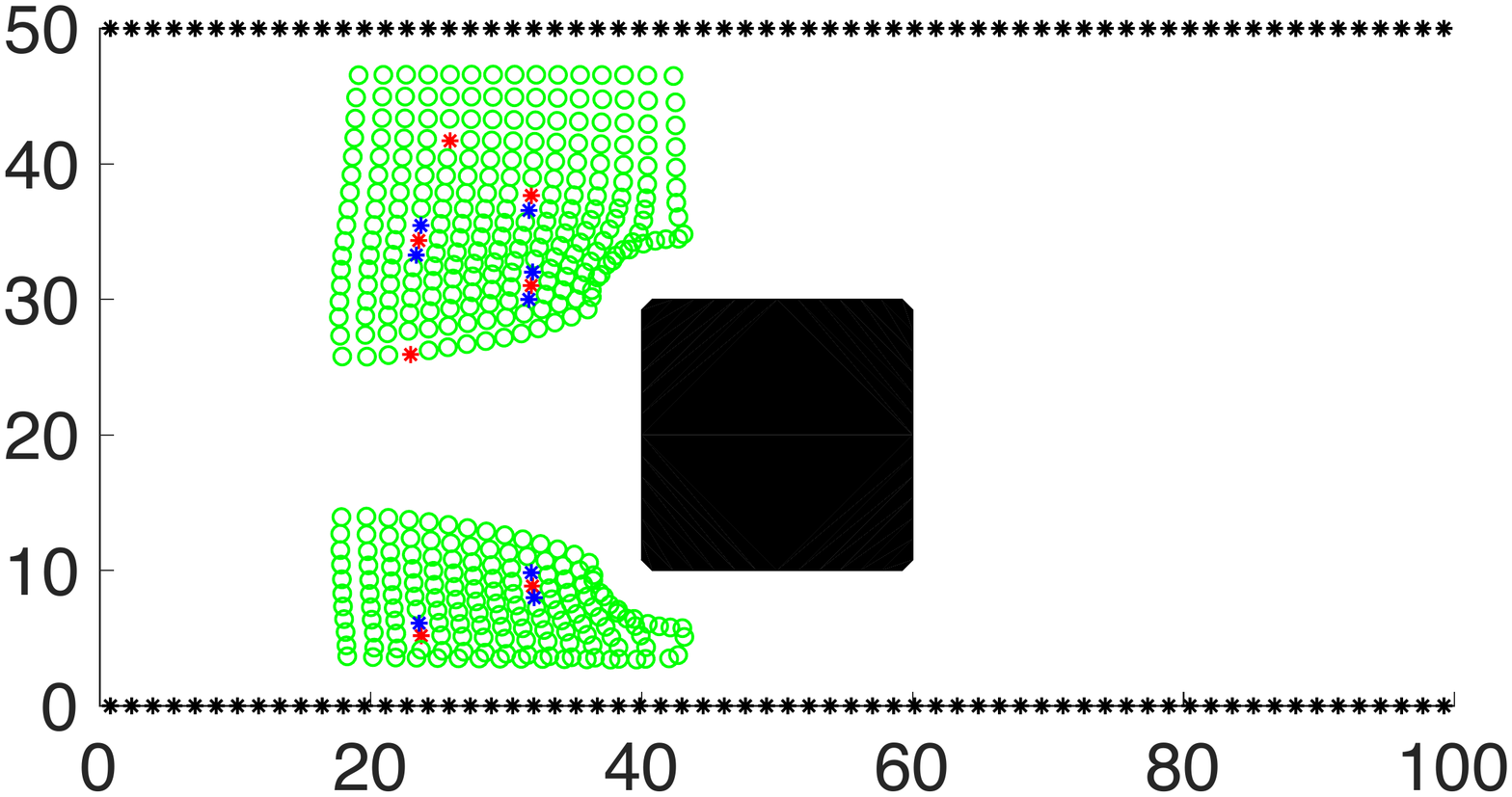}\\
		%	\vspace{-4cm}	
		\includegraphics[keepaspectratio=true, angle=90, width=0.45\textwidth]{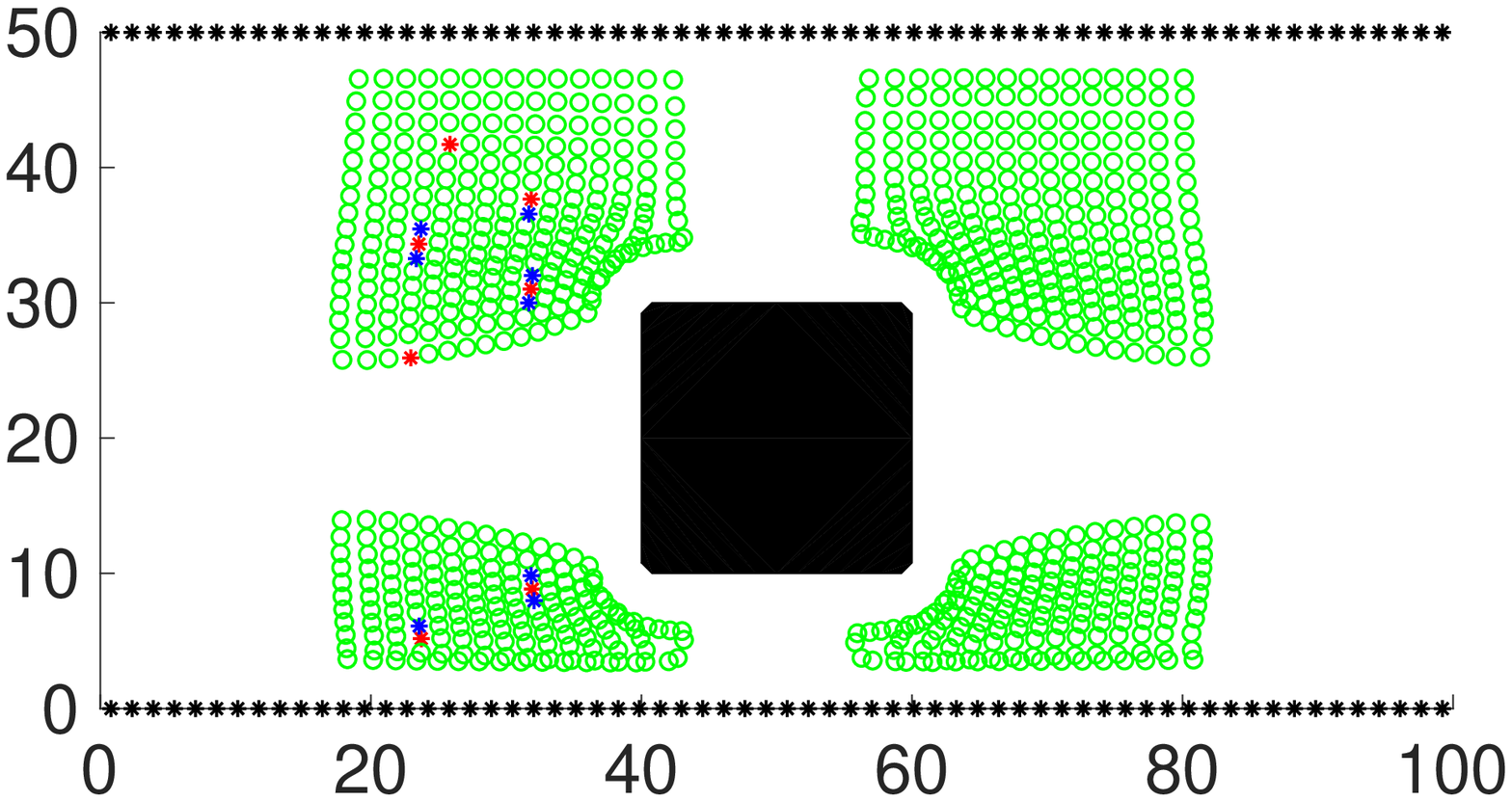}	 
		\includegraphics[keepaspectratio=true, angle=90, width=0.45\textwidth]{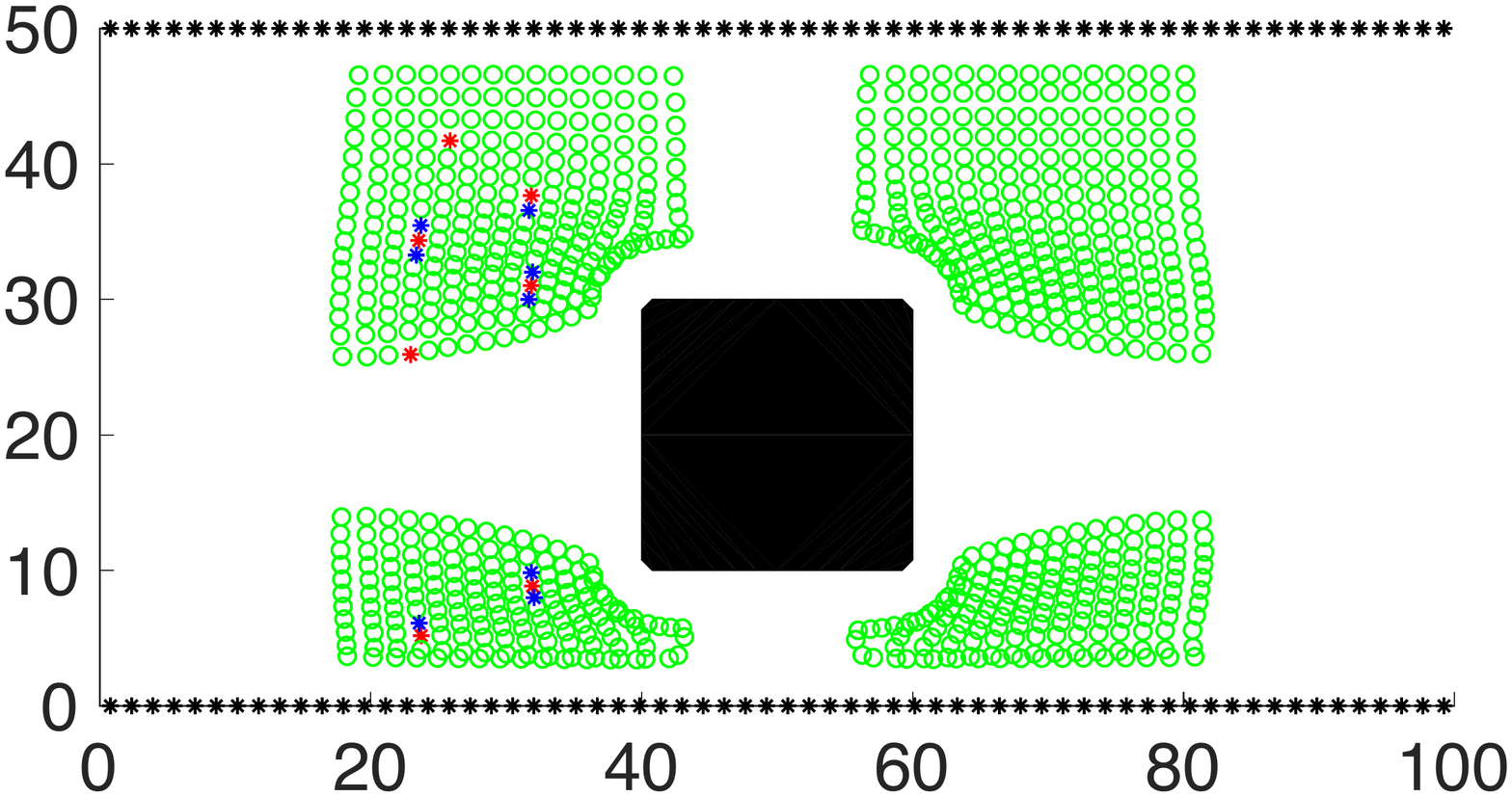}	 
		%	\vspace{-2cm}
		\caption{ Pedestrian dynamics  at time $t=10$ with influence (left) and without  influence (right) of contact time.   Row 1: Bi-directional flow. Row 2:  Uni-directional flow around obstacle. Row 3: Bi-directional flow around obstacle. Red indicate infected, green indicate susceptibles and blue indicate probably exposed pedestrians. }
		\label{phi_v_1_t10}
		\centering
	\end{figure}

\begin{figure}
		\centering
		\includegraphics[keepaspectratio=true, angle=90, width=0.45\textwidth]{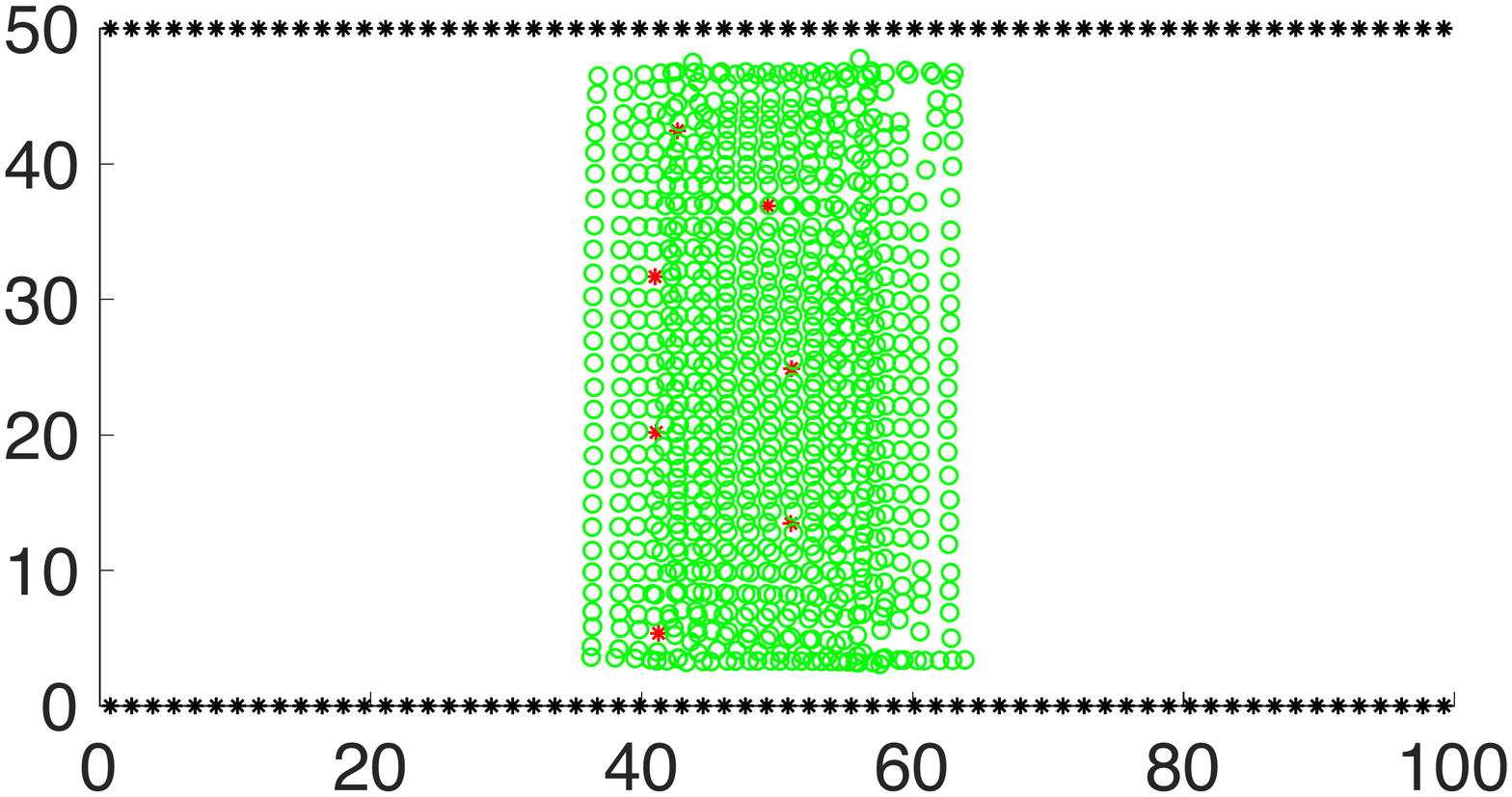}
		\includegraphics[keepaspectratio=true, angle=90, width=0.45\textwidth]{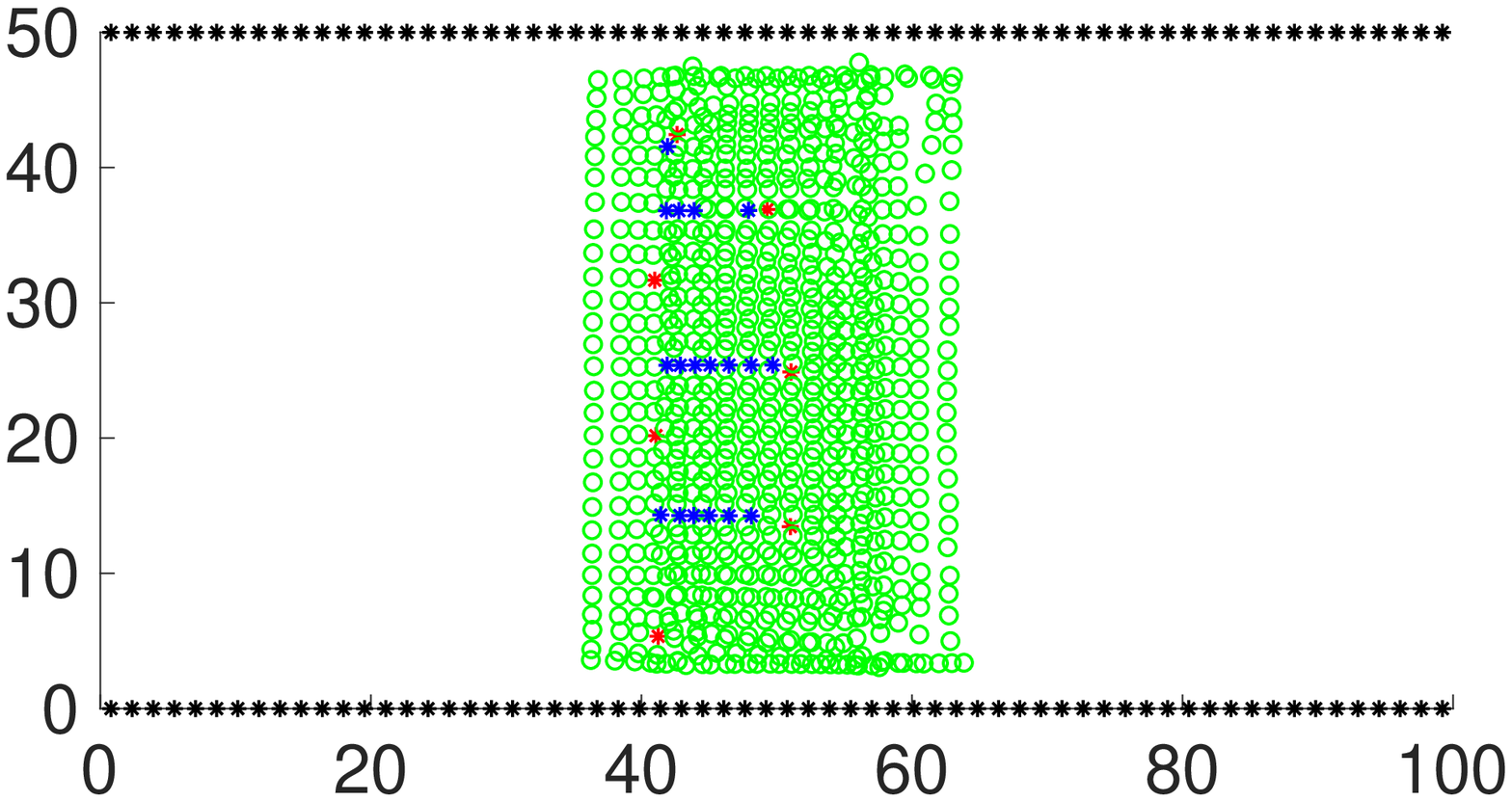}\\
		%%	\vspace{-4cm}	
		\includegraphics[keepaspectratio=true, angle=90, width=0.45\textwidth]{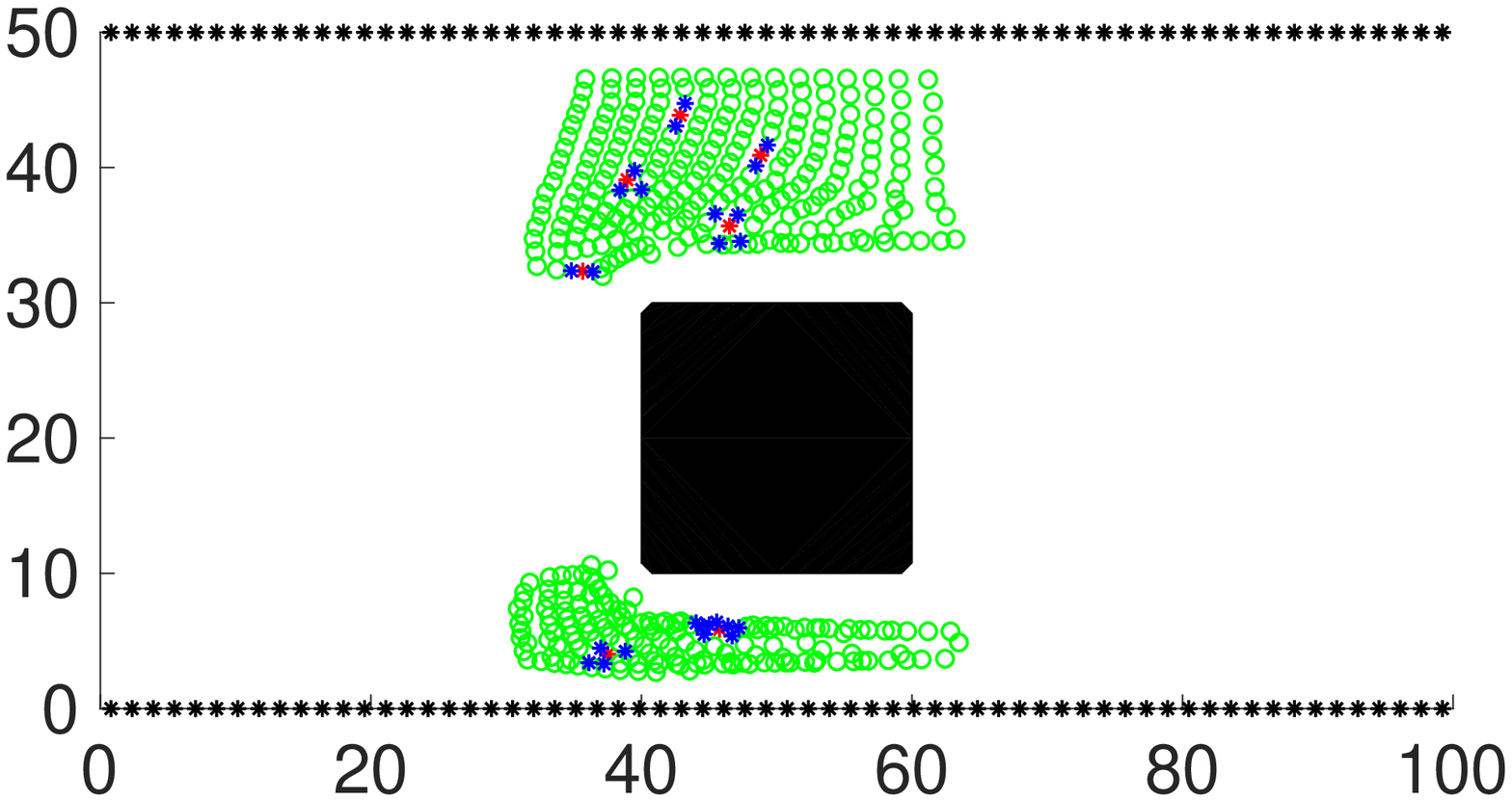}
		\includegraphics[keepaspectratio=true, angle=90, width=0.45\textwidth]{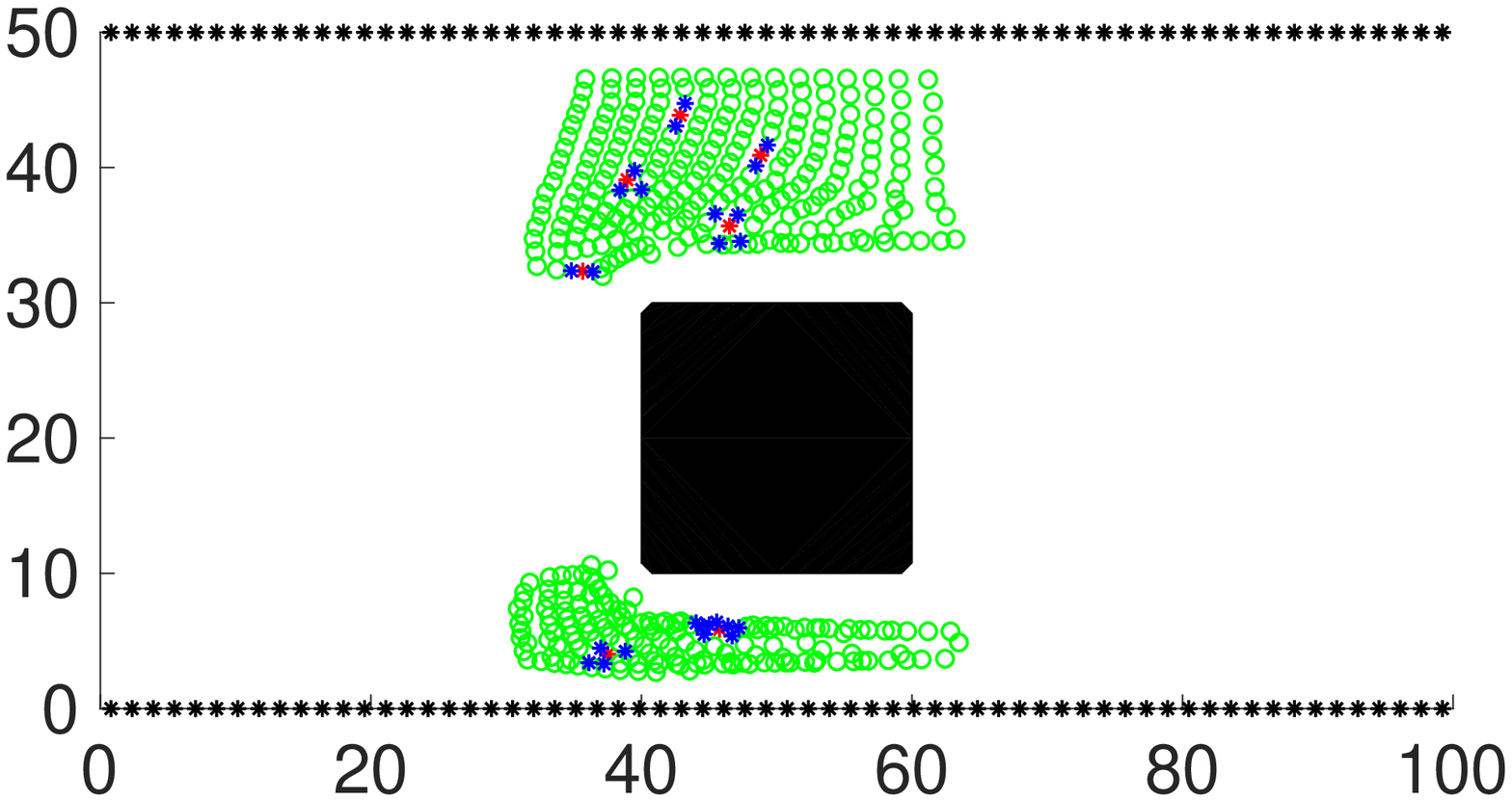}\\
		%%\vspace{-4cm}	
		\includegraphics[keepaspectratio=true, angle=90, width=0.45\textwidth]{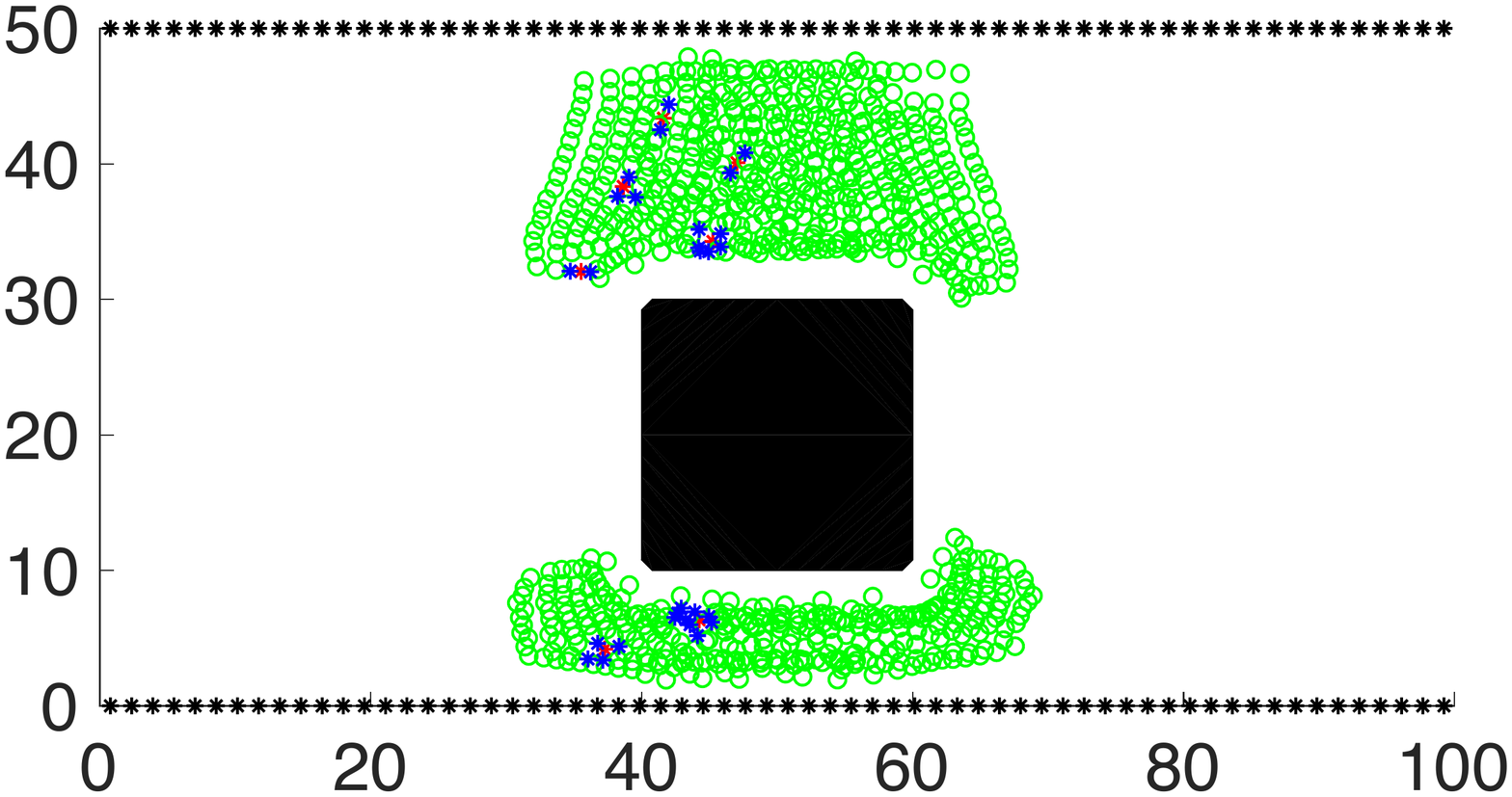}	 
		\includegraphics[keepaspectratio=true, angle=90, width=0.45\textwidth]{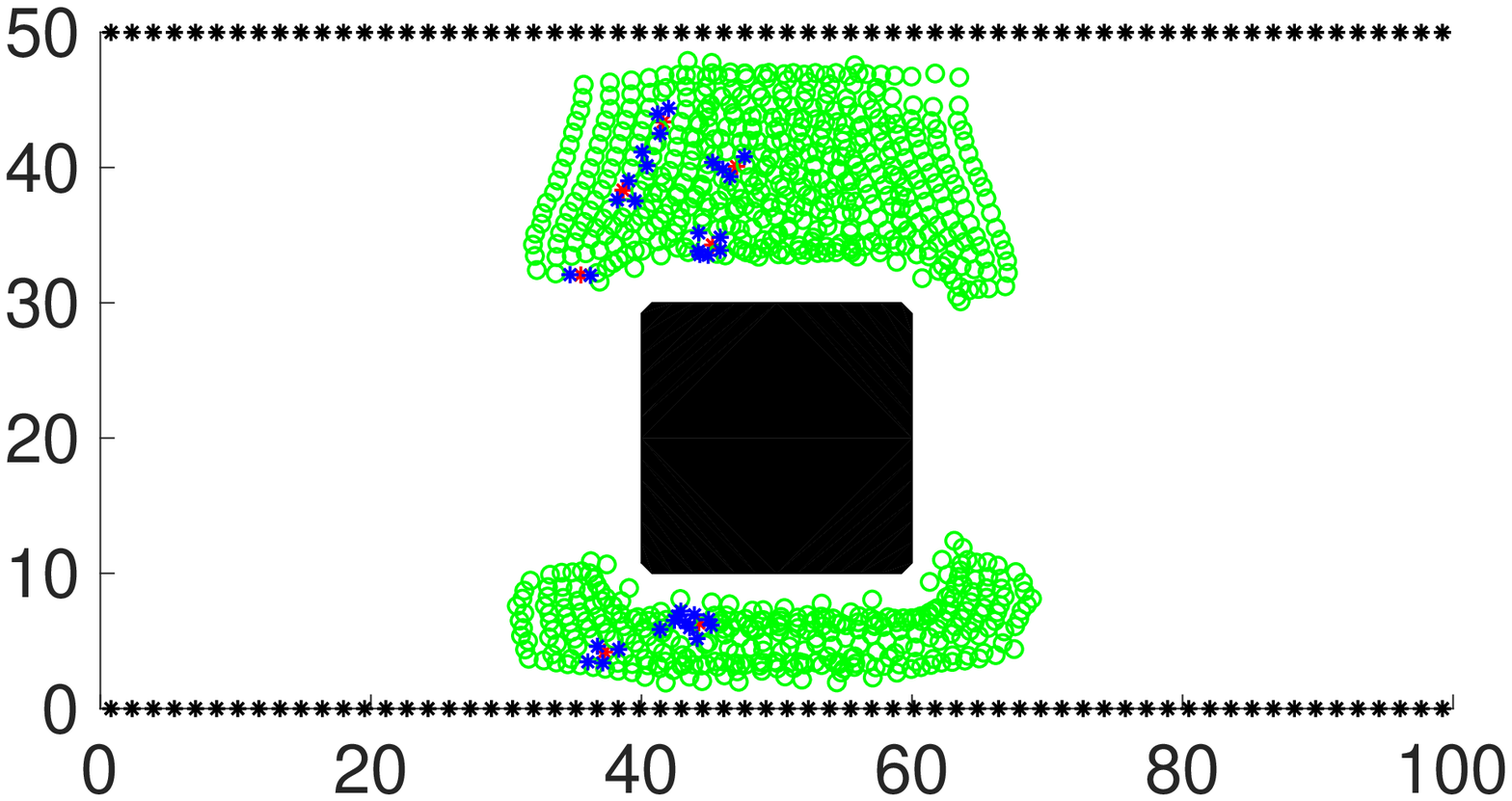}	 
		%%	\vspace{-2cm}
		\caption{ Pedestrian dynamics  at time $t=20$ with  influence (left) and without  influence (right) of contact time.   Row 1: Bi-directional flow. Row 2:  Uni-directional flow around obstacle. Row 3: Bi-directional flow around obstacle. Red indicate infected, green indicate susceptibles and blue indicate probably exposed pedestrians. }
		\label{phi_v_1_t20}
		\centering
	\end{figure}

\begin{figure}
		\centering
		\includegraphics[keepaspectratio=true, angle=90, width=0.45\textwidth]{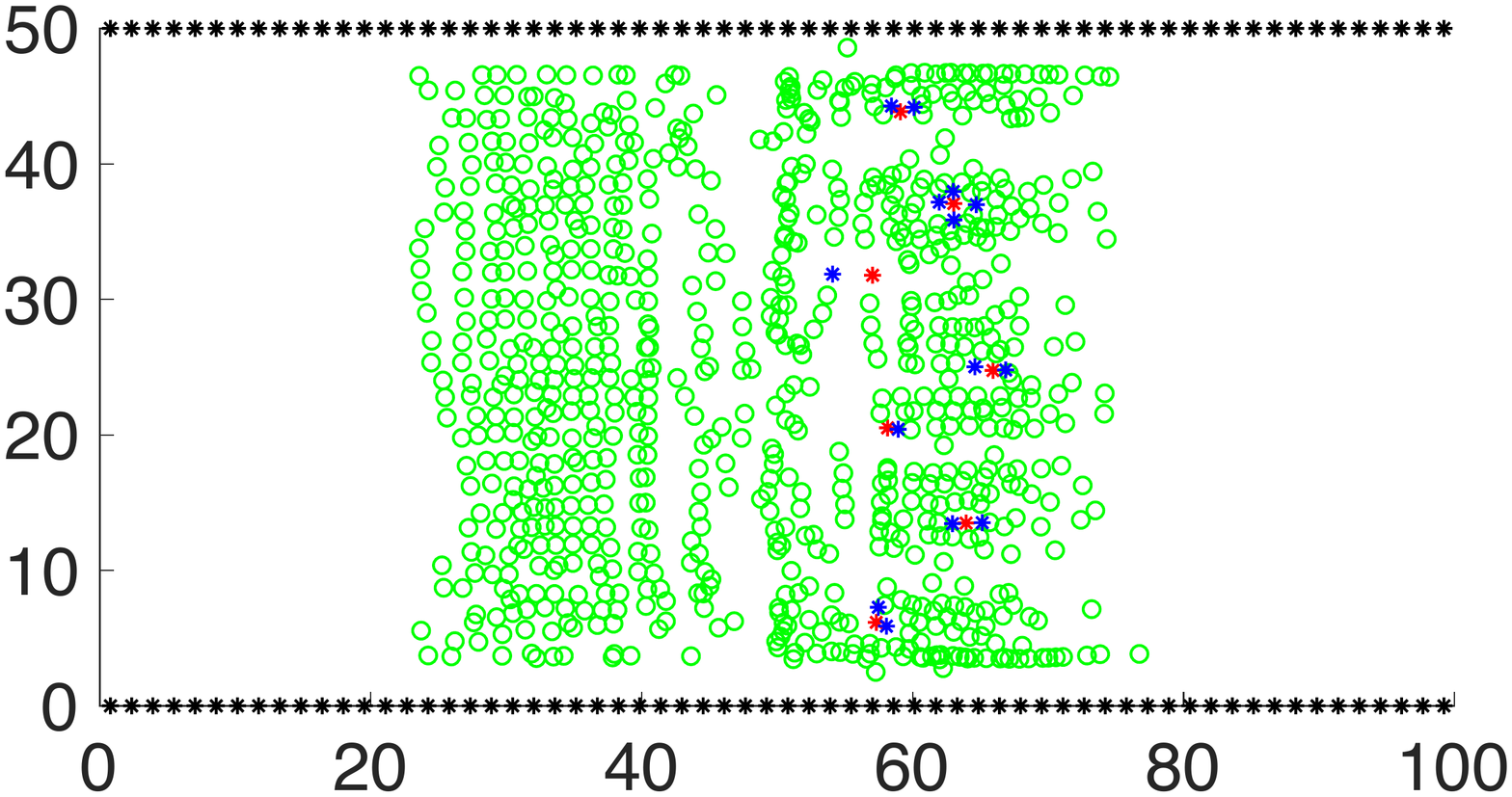}
		\includegraphics[keepaspectratio=true, angle=90, width=0.45\textwidth]{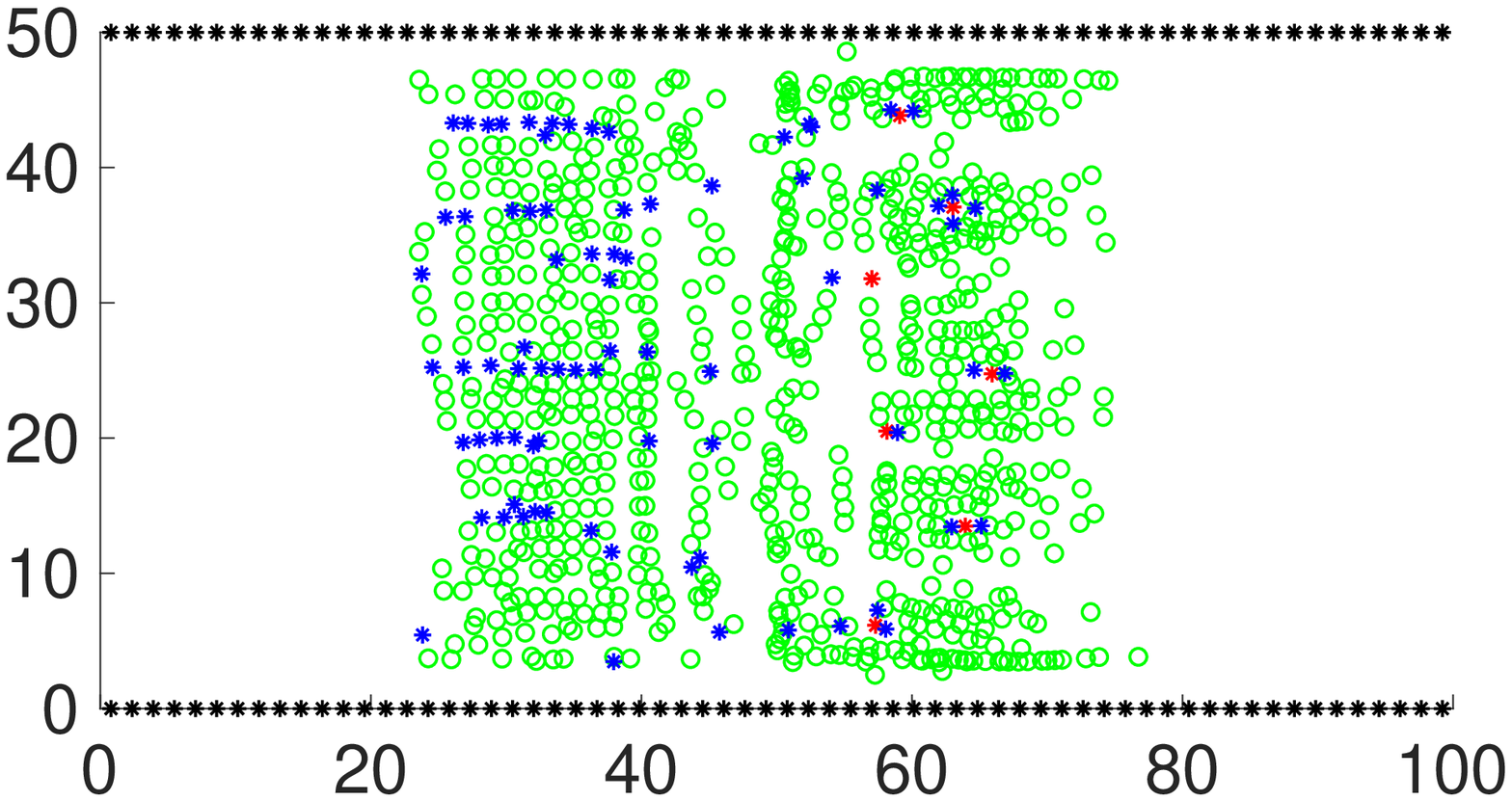}\\
		%%	\vspace{-4cm}	
		\includegraphics[keepaspectratio=true, angle=90, width=0.45\textwidth]{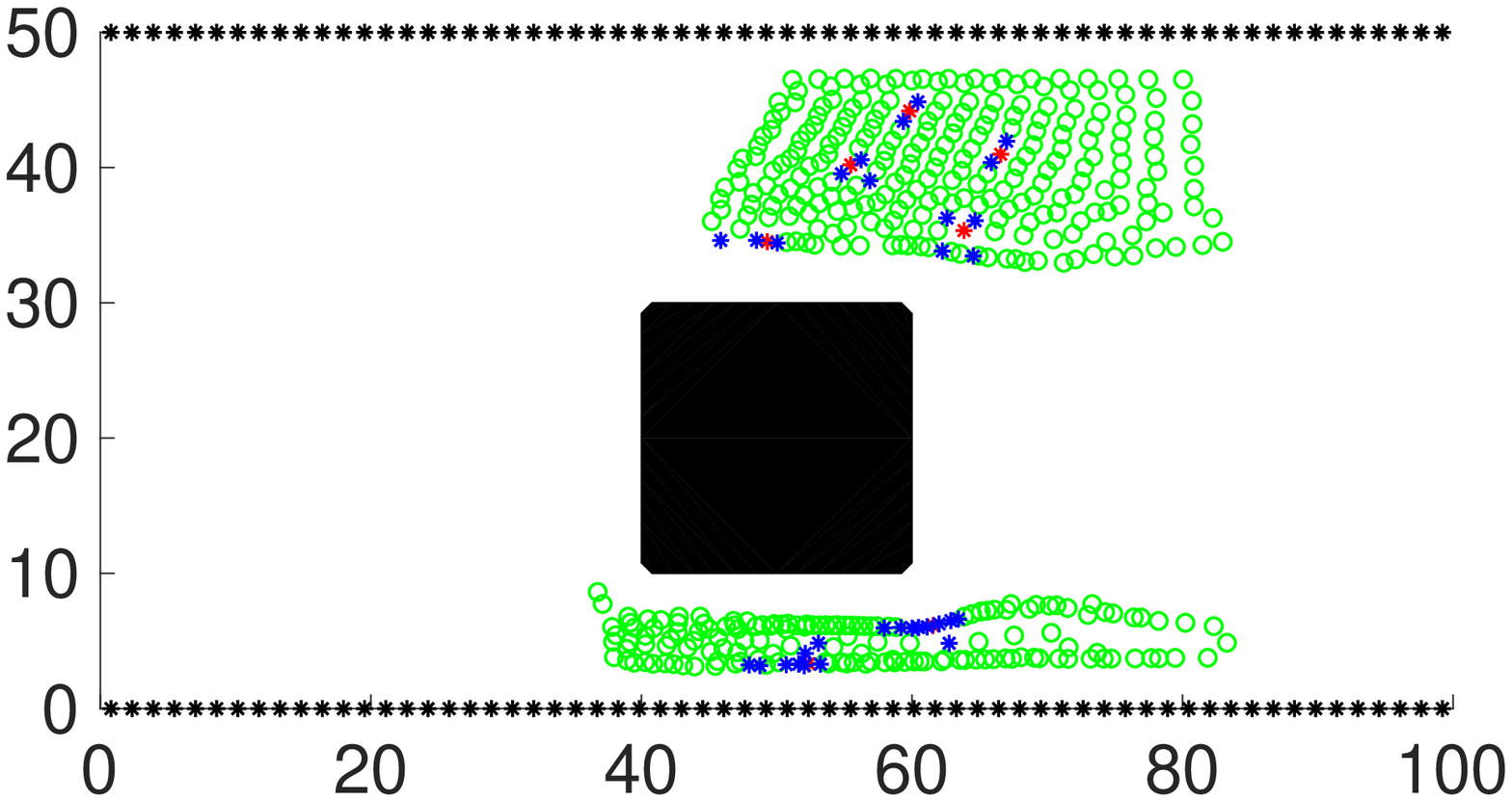}
		\includegraphics[keepaspectratio=true, angle=90, width=0.45\textwidth]{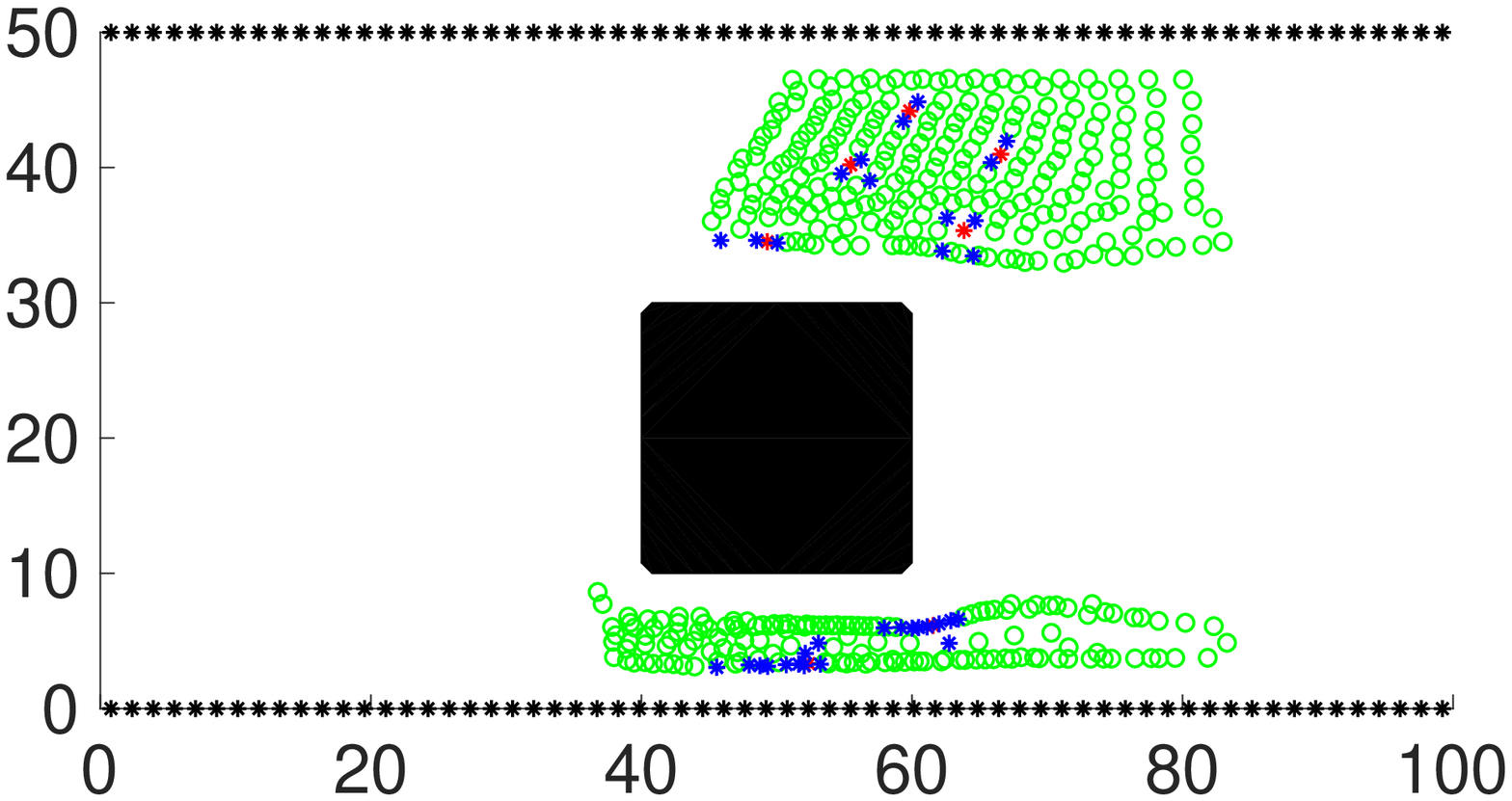}\\
		%%	\vspace{-4cm}	
		\includegraphics[keepaspectratio=true, angle=90, width=0.45\textwidth]{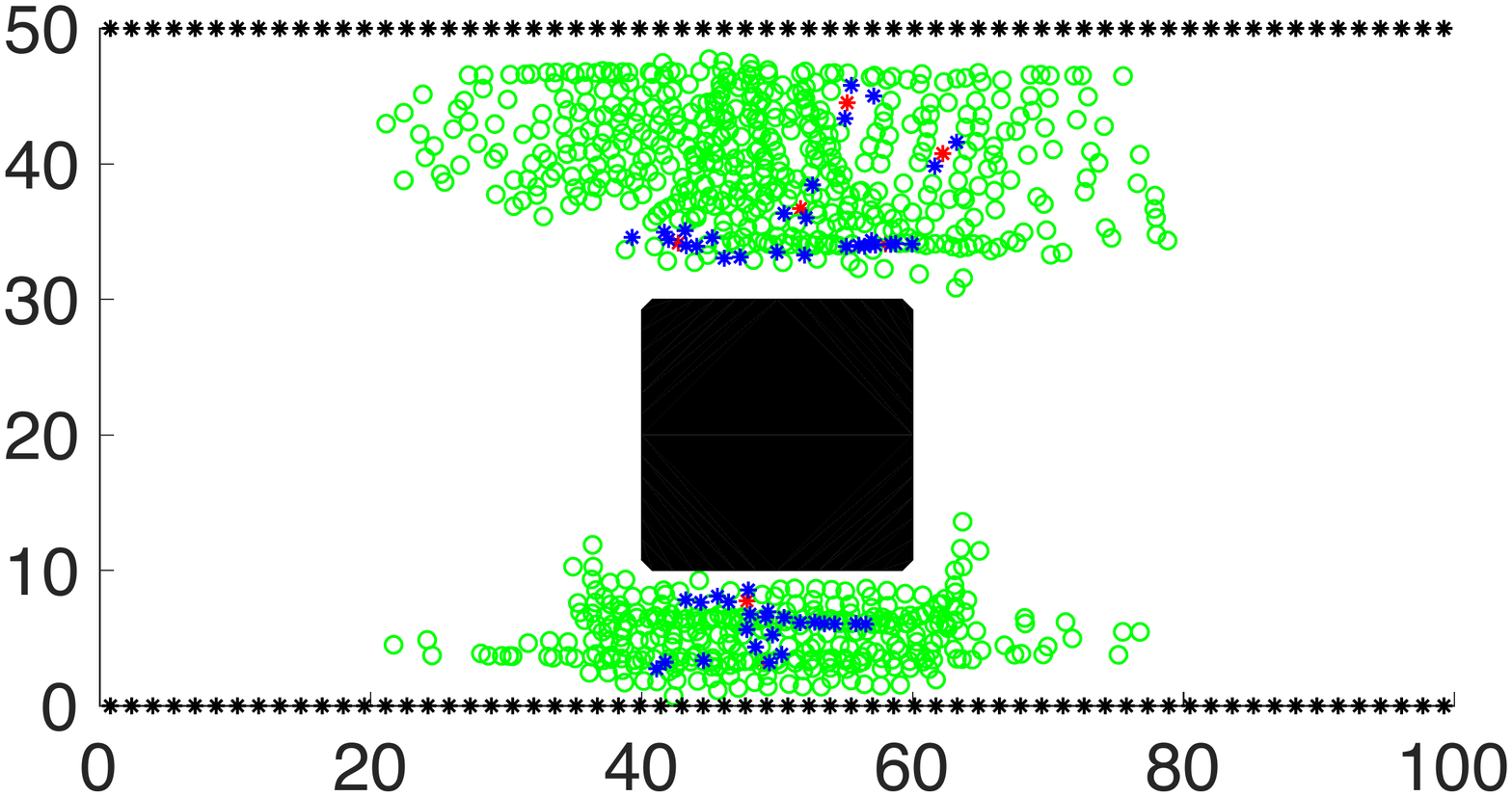}	 
		\includegraphics[keepaspectratio=true, angle=90, width=0.45\textwidth]{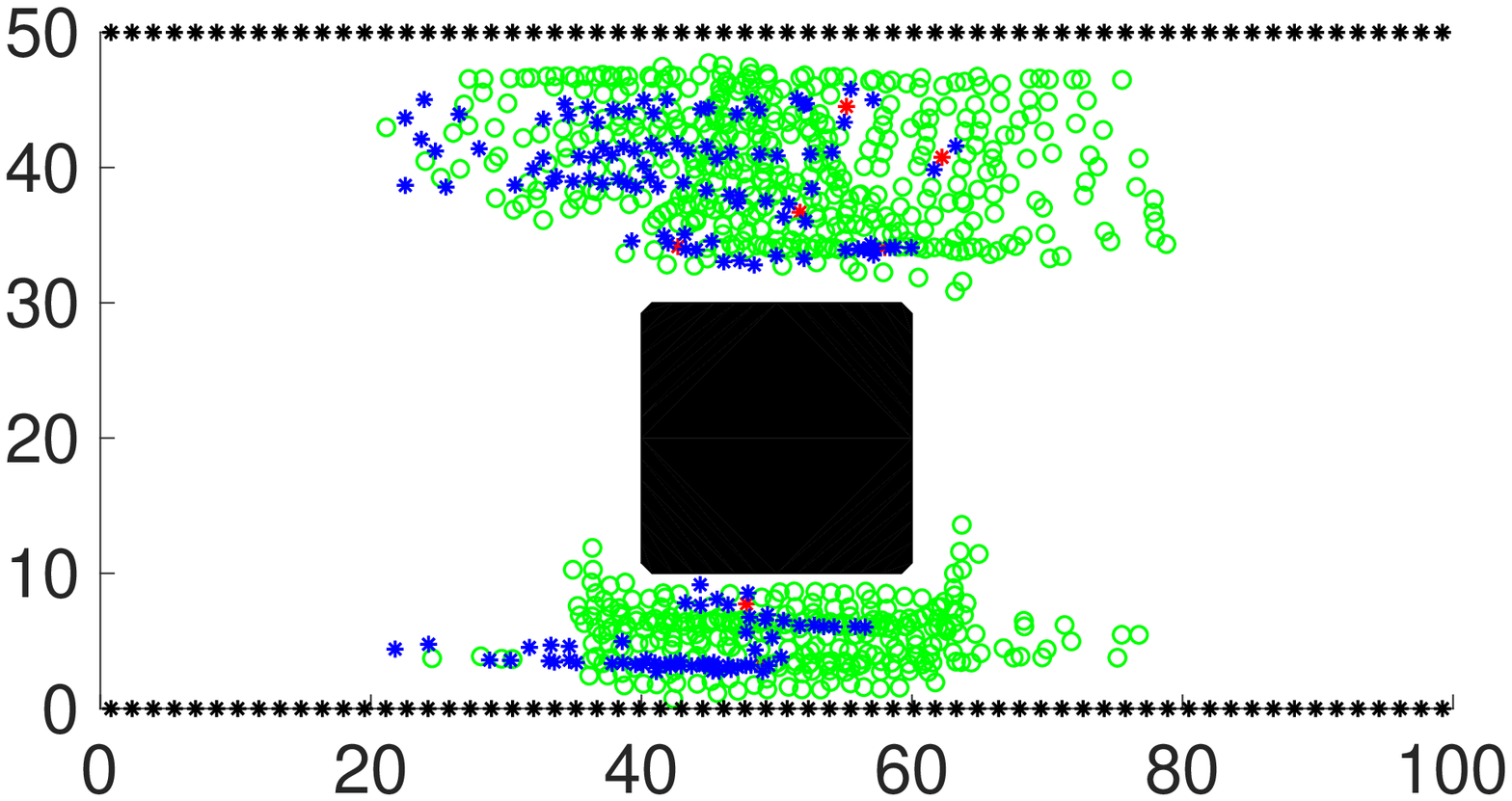}	 
		%%	\vspace{-2cm}
		\caption{ Pedestrian dynamics  at time $t=30$ with  influence (left) and without  influence (right) of contact time.   Row 1: Bi-directional flow. Row 2:  Uni-directional flow around obstacle. Row 3: Bi-directional flow around obstacle. Red indicate infected, green indicate susceptibles and blue indicate probably exposed pedestrians. }
		\label{phi_v_1_t30}
		\centering
	\end{figure}

\begin{figure}
		\centering
		\includegraphics[keepaspectratio=true, angle=90, width=0.45\textwidth]{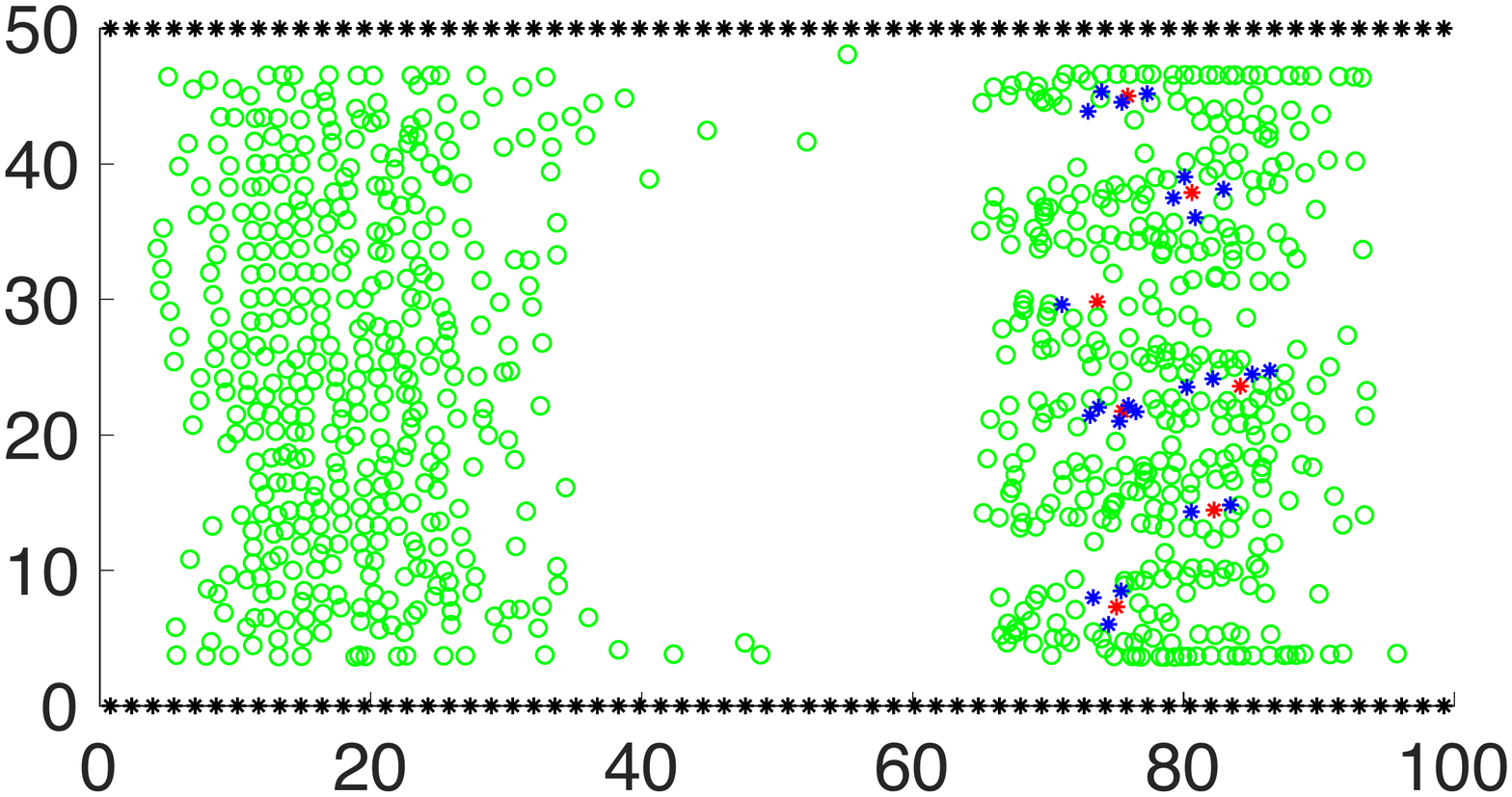}
		\includegraphics[keepaspectratio=true, angle=90, width=0.45\textwidth]{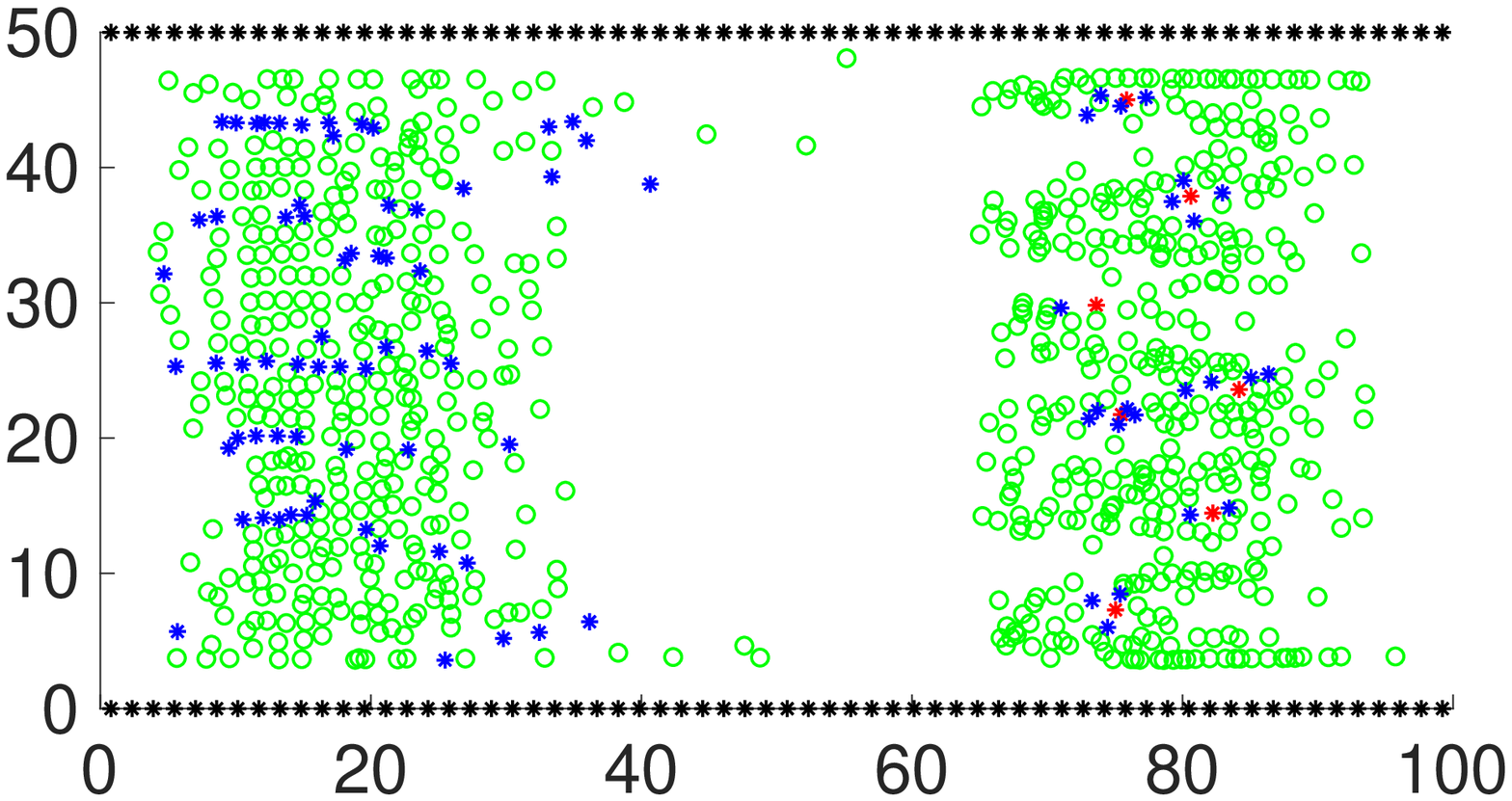}\\
		%%	\vspace{-4cm}	
		\includegraphics[keepaspectratio=true, angle=90, width=0.45\textwidth]{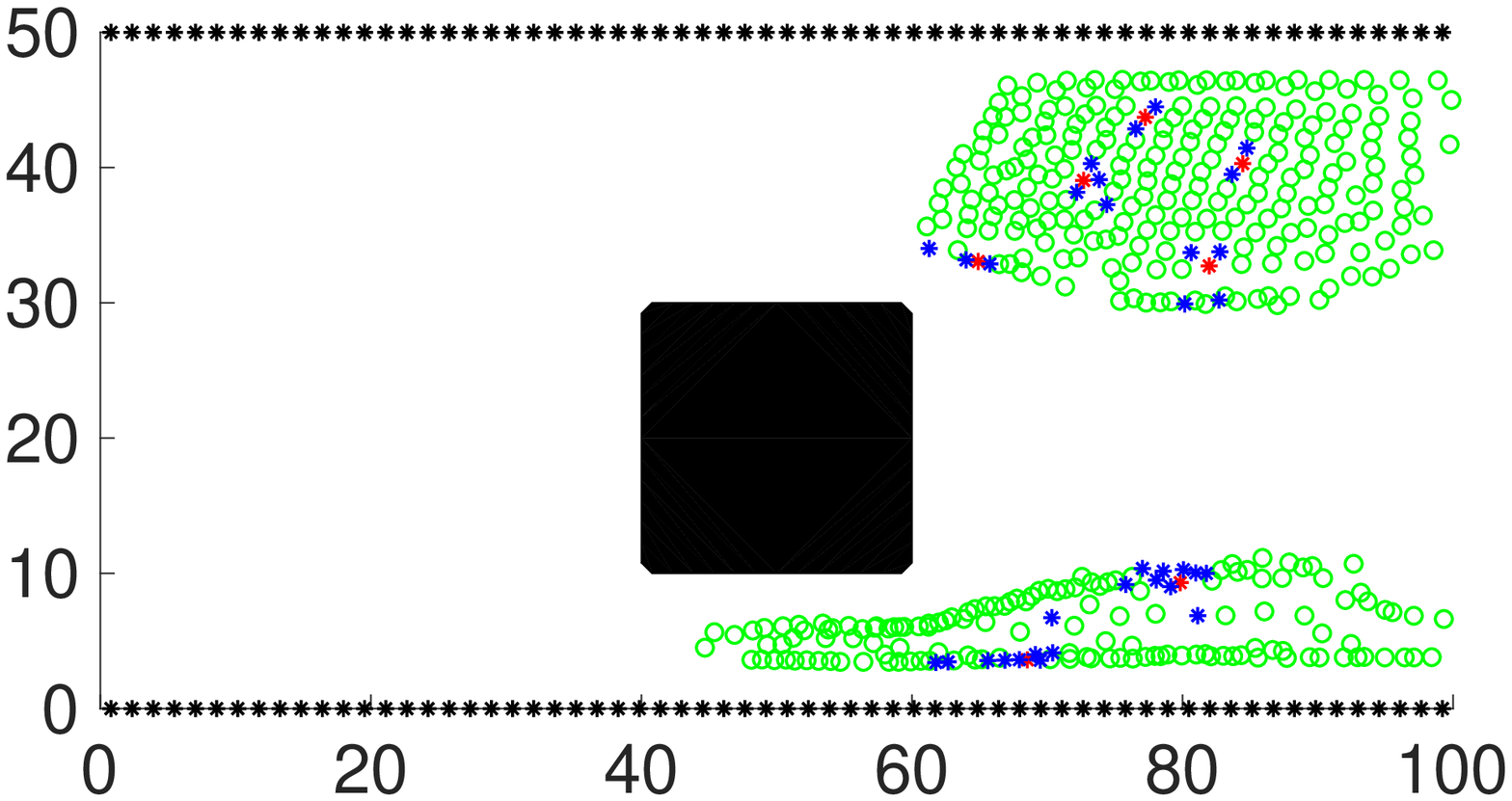}
		\includegraphics[keepaspectratio=true, angle=90, width=0.45\textwidth]{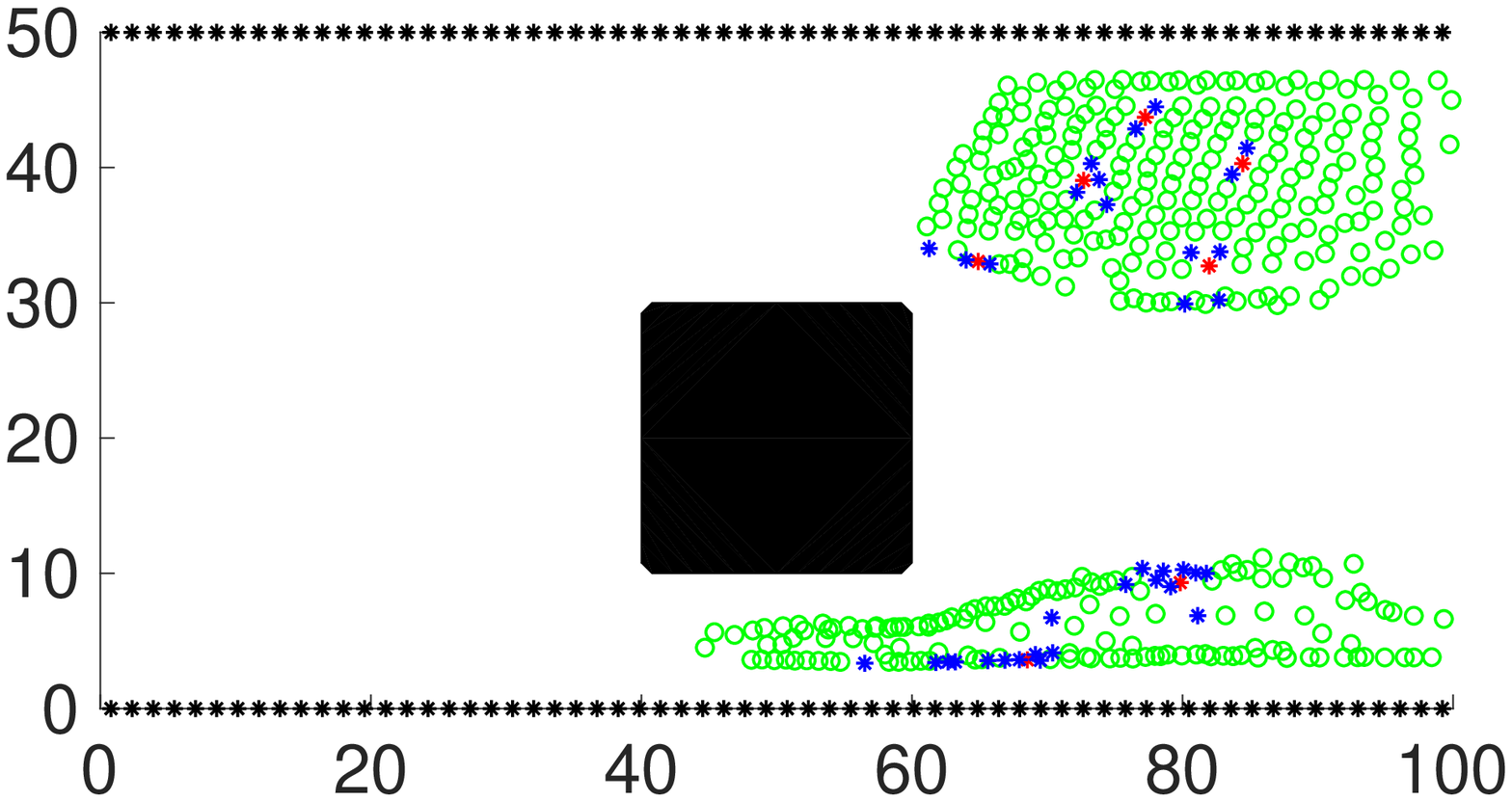}\\
		%%	\vspace{-4cm}	
		\includegraphics[keepaspectratio=true, angle=90, width=0.45\textwidth]{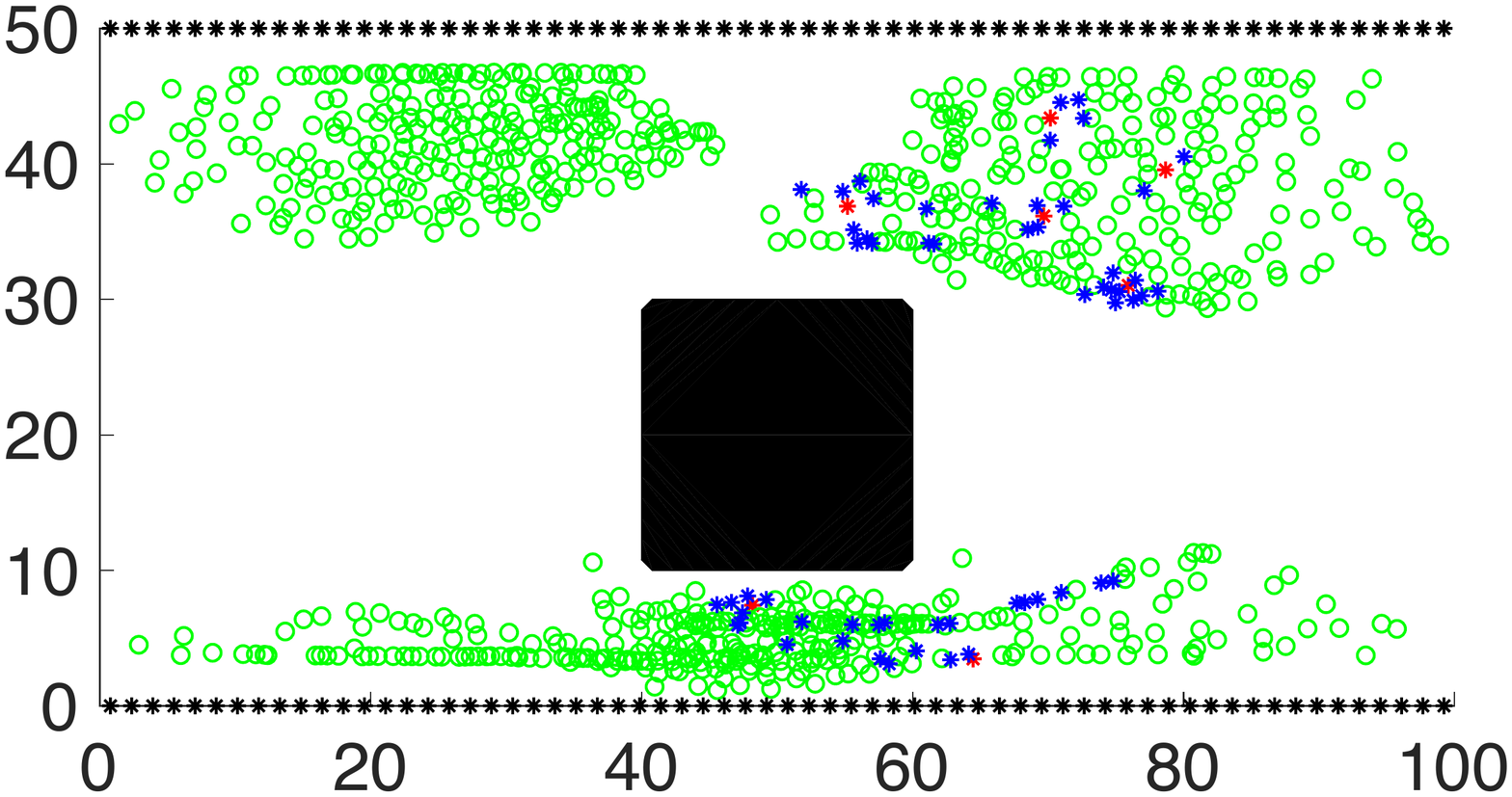}	 
		\includegraphics[keepaspectratio=true, angle=90, width=0.45\textwidth]{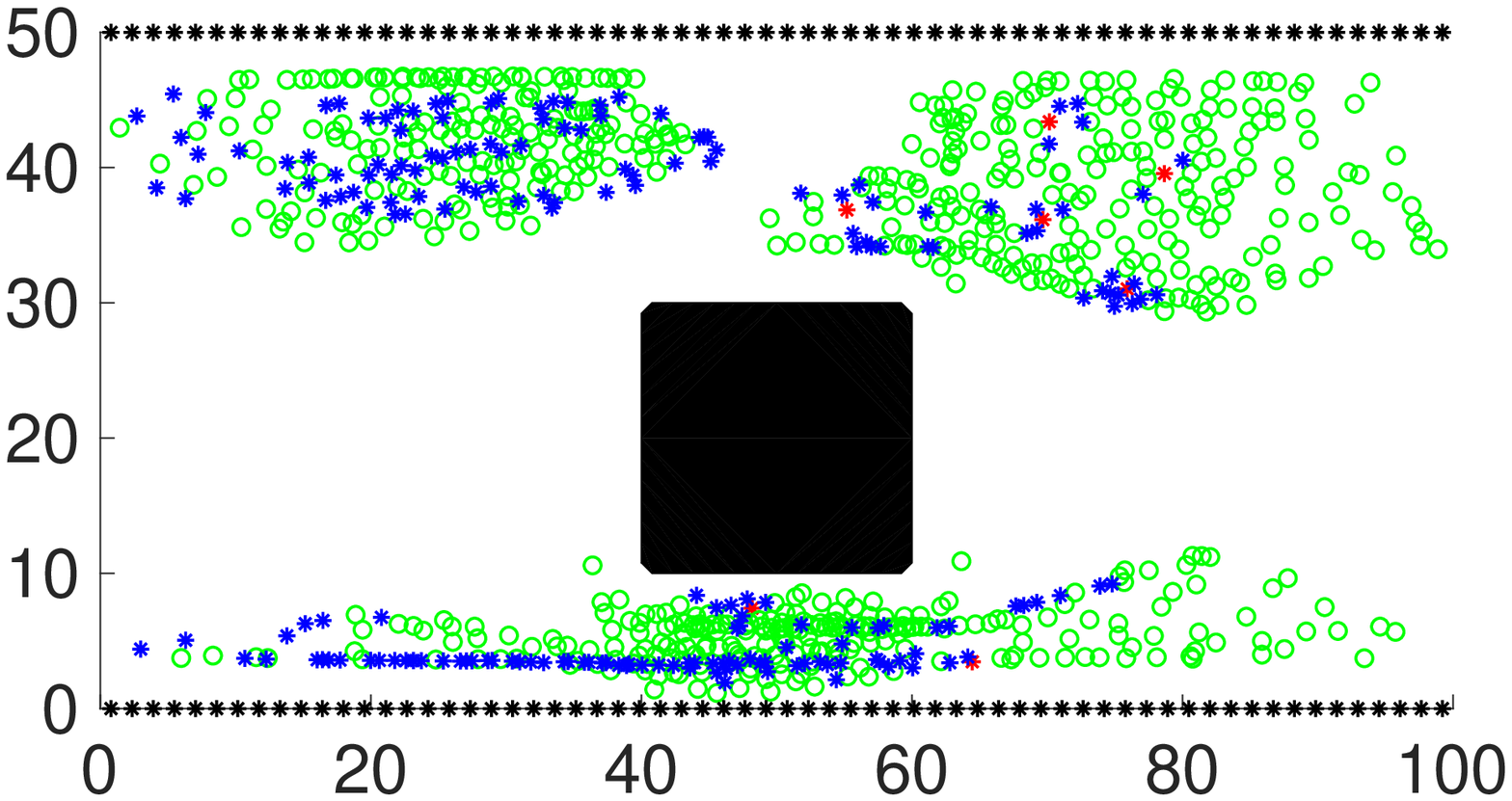}	 
		%%	\vspace{-2cm}
		\caption{ Pedestrian dynamics  at time $t=40$ with  influence (left) and without  influence (right) of contact time.   Row 1: Bi-directional flow. Row 2:  Uni-directional flow around obstacle. Row 3: Bi-directional flow around obstacle. Red indicate infected, green indicate susceptibles and blue indicate probably exposed pedestrians. }
		\label{phi_v_1_t40}
		\centering
	\end{figure}

	Figures \ref{phi_v_1_t10} to \ref{phi_v_1_t40} 
	show the time evolution of  the moving grid points and the associated infection labels with influence (left column) and without  influence (right column) of contact time.  
	Row 1 shows bi-directional flow without an obstacle, row 2 shows uni-directional flow around an obstacle and row 3
	shows  bi-directional flow around  an obstacle. Red indicate infected, green indicate susceptibles and blue indicate exposed pedestrians.
	
	Here one observes, e.g. in Figure \ref{phi_v_1_t30}  or \ref{phi_v_1_t40}, that in  situations with bi-directional flow (top- and bottom row), the number of exposed patients is strongly reduced, if the contact time is taken  into account.
	For uni-directional flow, the differences are much smaller as expected, since pedestrians stay near to each other for a longer time during the evolution.
	
	Comparing row 2 (uni-dirctional with obstacle) with the uni-directional case without obstacle (no exposed pedestrians), one observes that the number of exposed patients  is considerably increased due to the denser pedestrian crowd surrounding the obstacle.
	Similar observations can be made  comparing row 1 and row 3.

	%%%%%%%%%%
	
	%%%%%%%%%%%
	
\begin{figure}
		\centering
		\includegraphics[keepaspectratio=true, angle=90, width=0.45\textwidth]{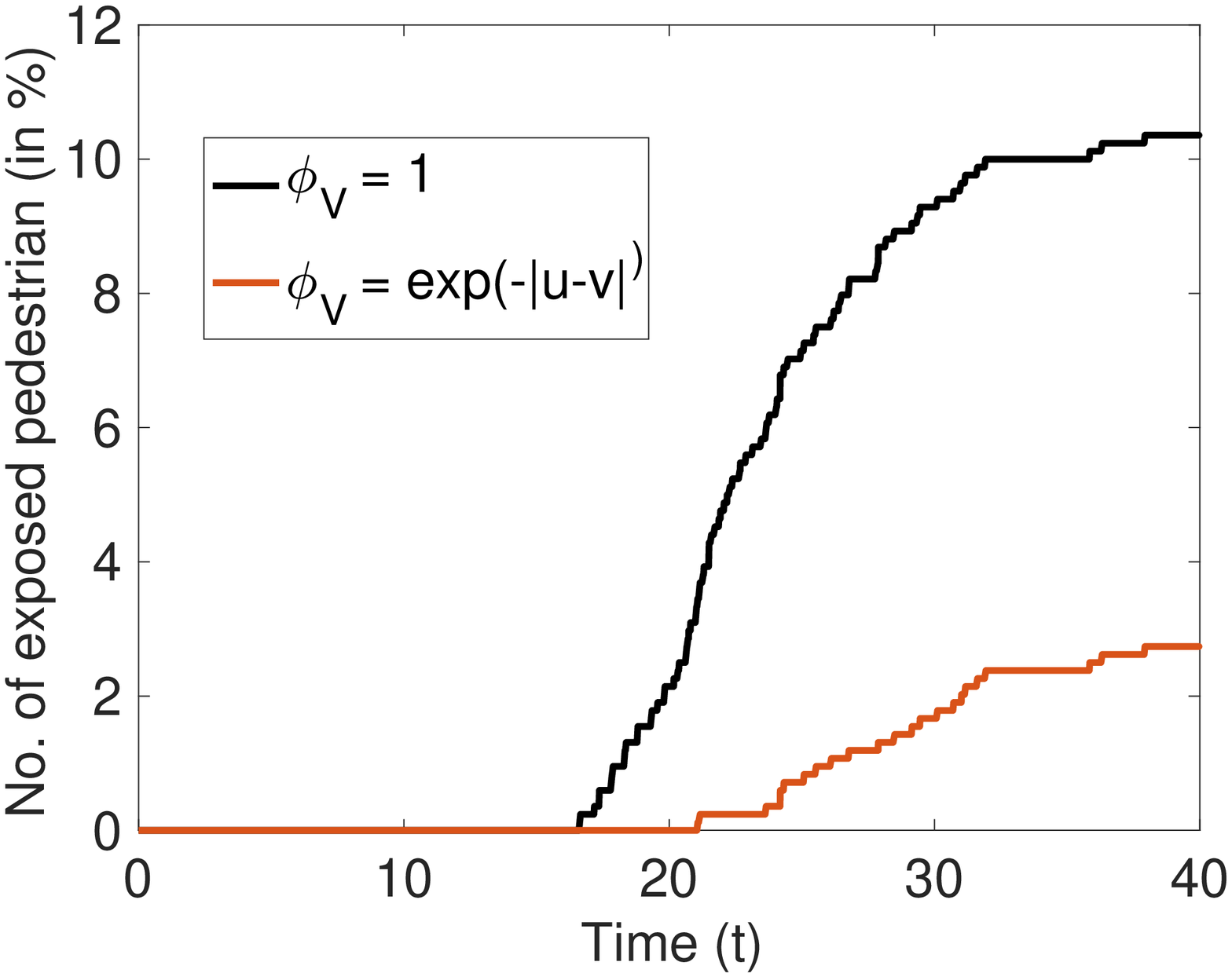}	
		\includegraphics[keepaspectratio=true, angle=90, width=0.45\textwidth]{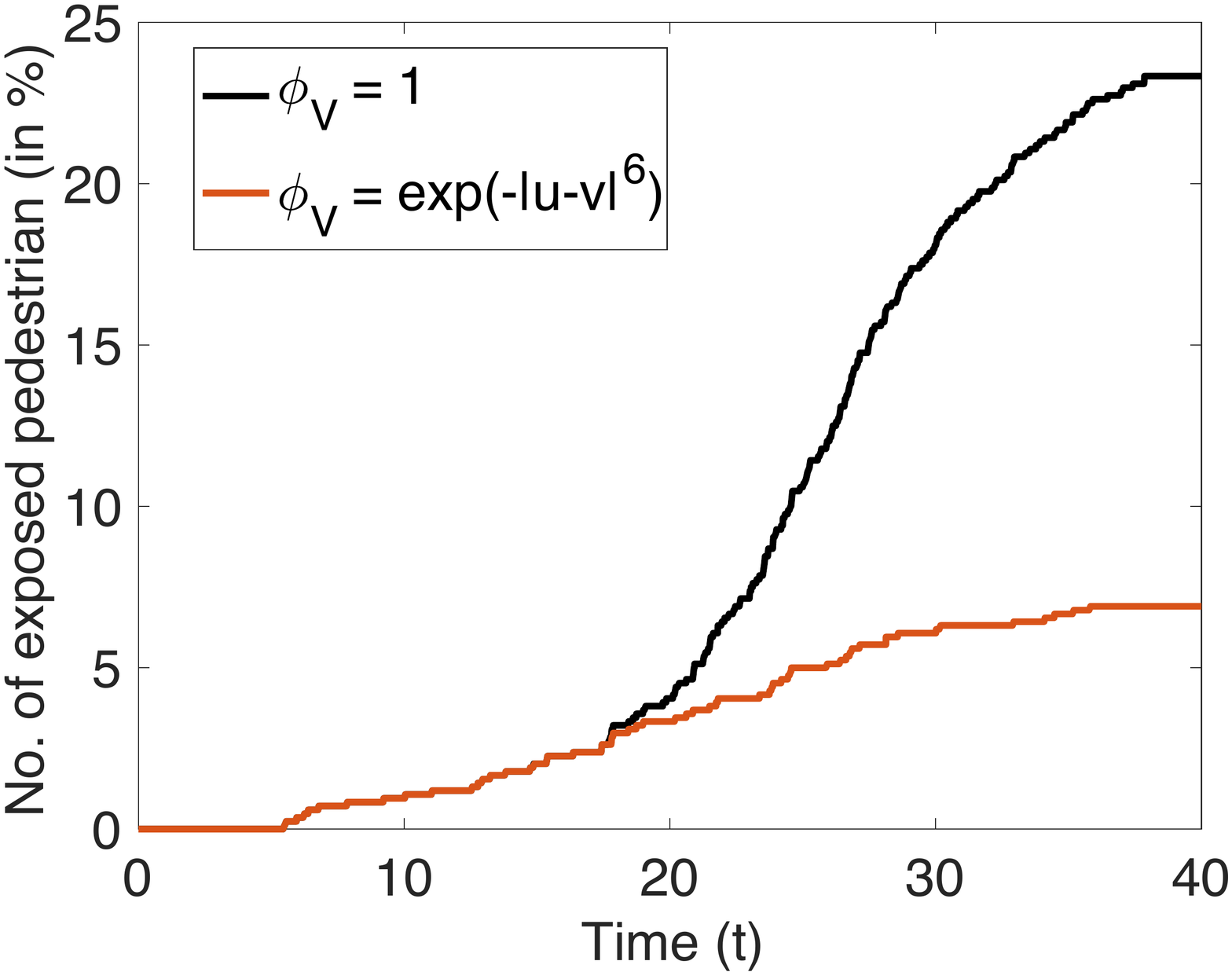}	
		\includegraphics[keepaspectratio=true, angle=90, width=0.45\textwidth]{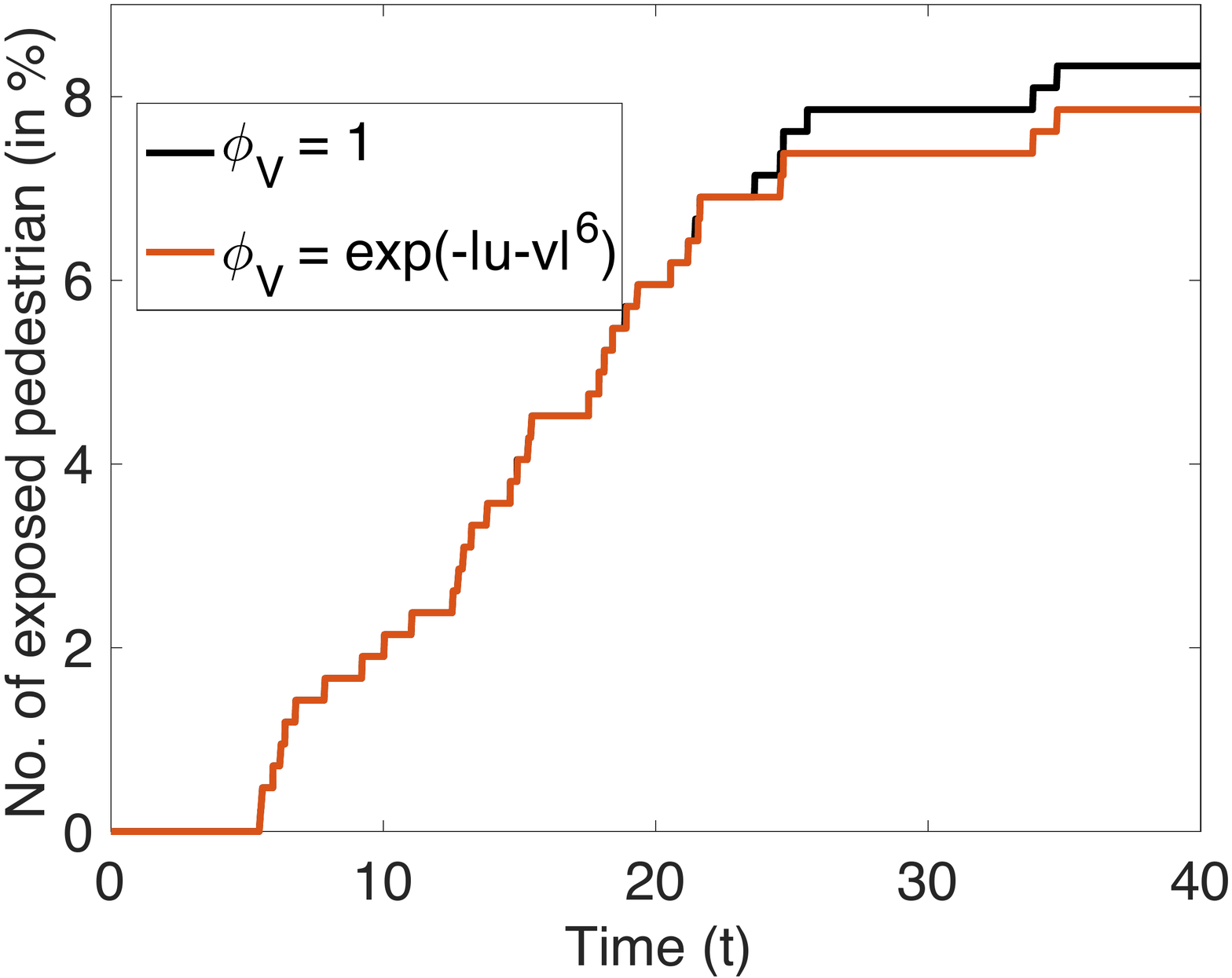}	
		\caption{Number of  pedestrian with an increased  probablity of being exposed (in \%) vs time. First row bi-directional flow without obstacle (left) and with obstacle (right). Second row: uni-directional with obstacle. }
		\label{exposed_ped}
		\centering
	\end{figure}

	%%%%%%%%%%%%%%%%%
	In Figure \ref{exposed_ped} we have plotted the number of  pedestrians with a higher probability of being exposed versus time. One can observe, that the number of these  pedestrians is much higher if the influence of the contact time is neglected. This happens mainly in bi-directional flows, which is as expected, since, even  if pedestrians are coming close to each other, they pass each other  quickly and  the contact time is  short such that a contagion is less probable. On the other hand if pedestrians are walking in the same direction, the effect of neglecting the contact time is comparably  small. We mention, that the contagion model is based on very simple assumptions and obviously the parameters of the contagion model have to be adapted to experimental findings.
	We refer again to \cite{KQ} for similiar investigations.

	%%%%%%%%%%%%%%%%%%
	\begin{figure}
		\centering
		\includegraphics[keepaspectratio=true, angle=90, width=0.45\textwidth]{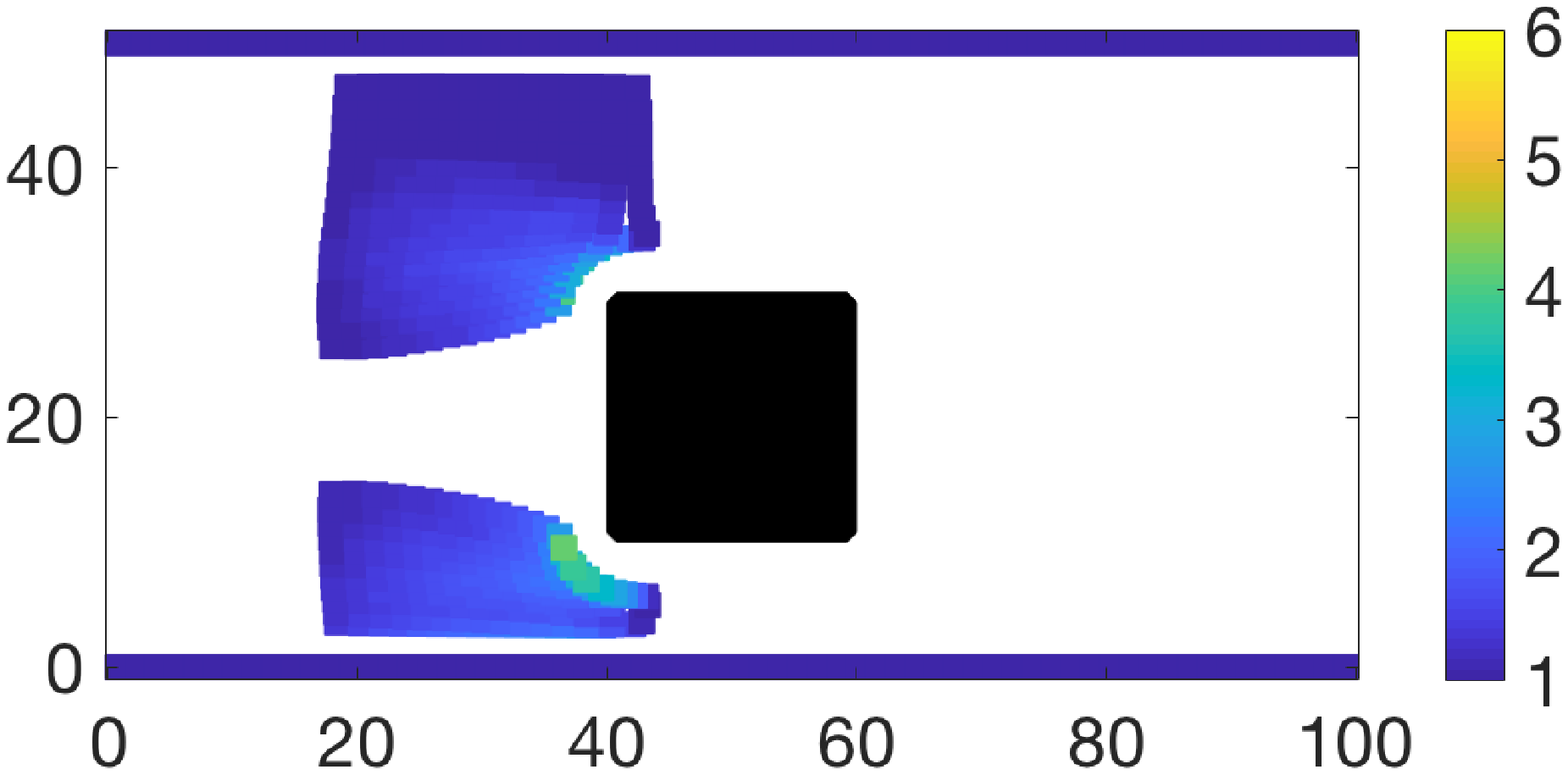}
		\includegraphics[keepaspectratio=true, angle=90, width=0.45\textwidth]{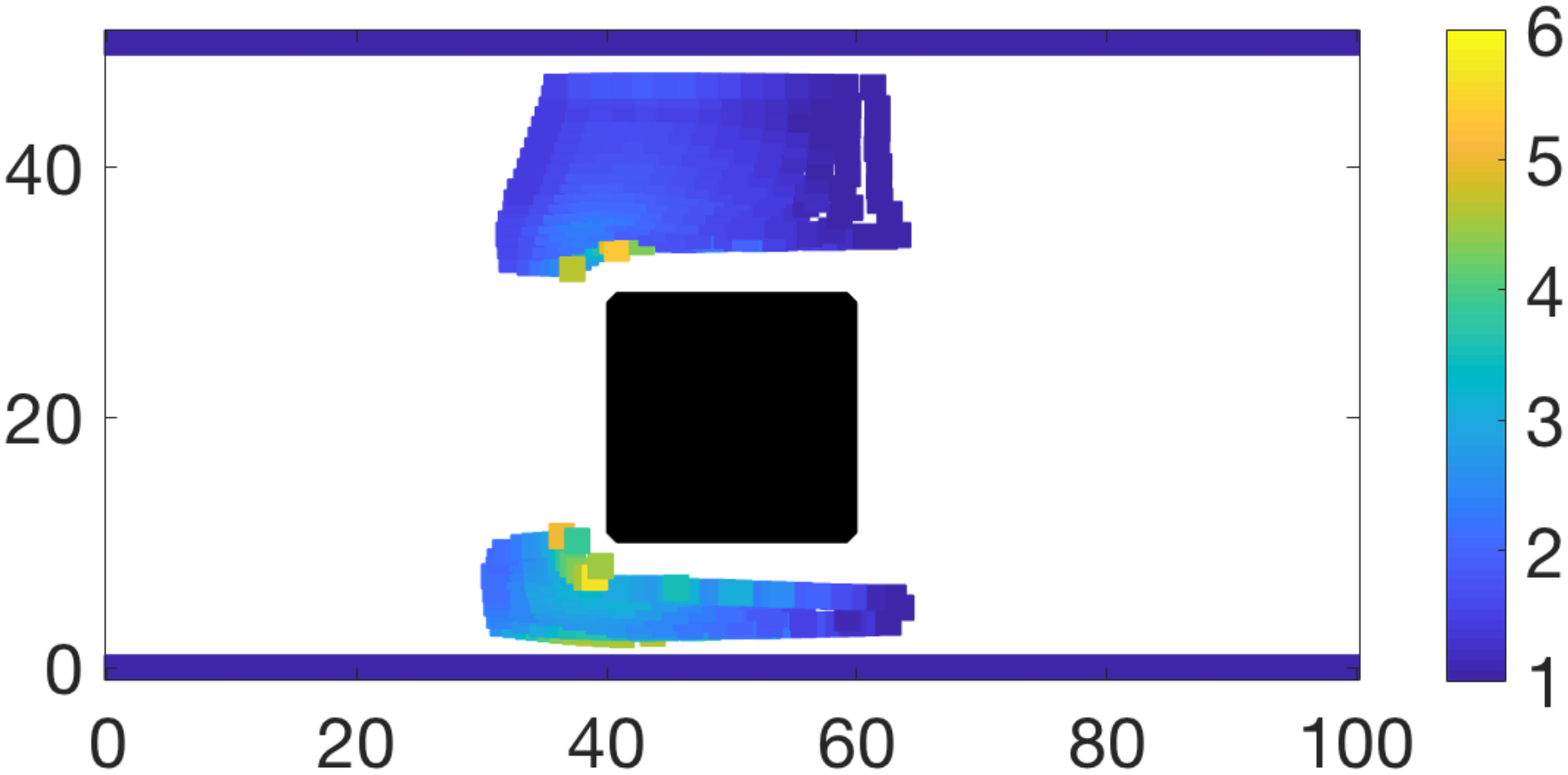}\\
		%	\vspace{-4cm}	
		\includegraphics[keepaspectratio=true, angle=90, width=0.45\textwidth]{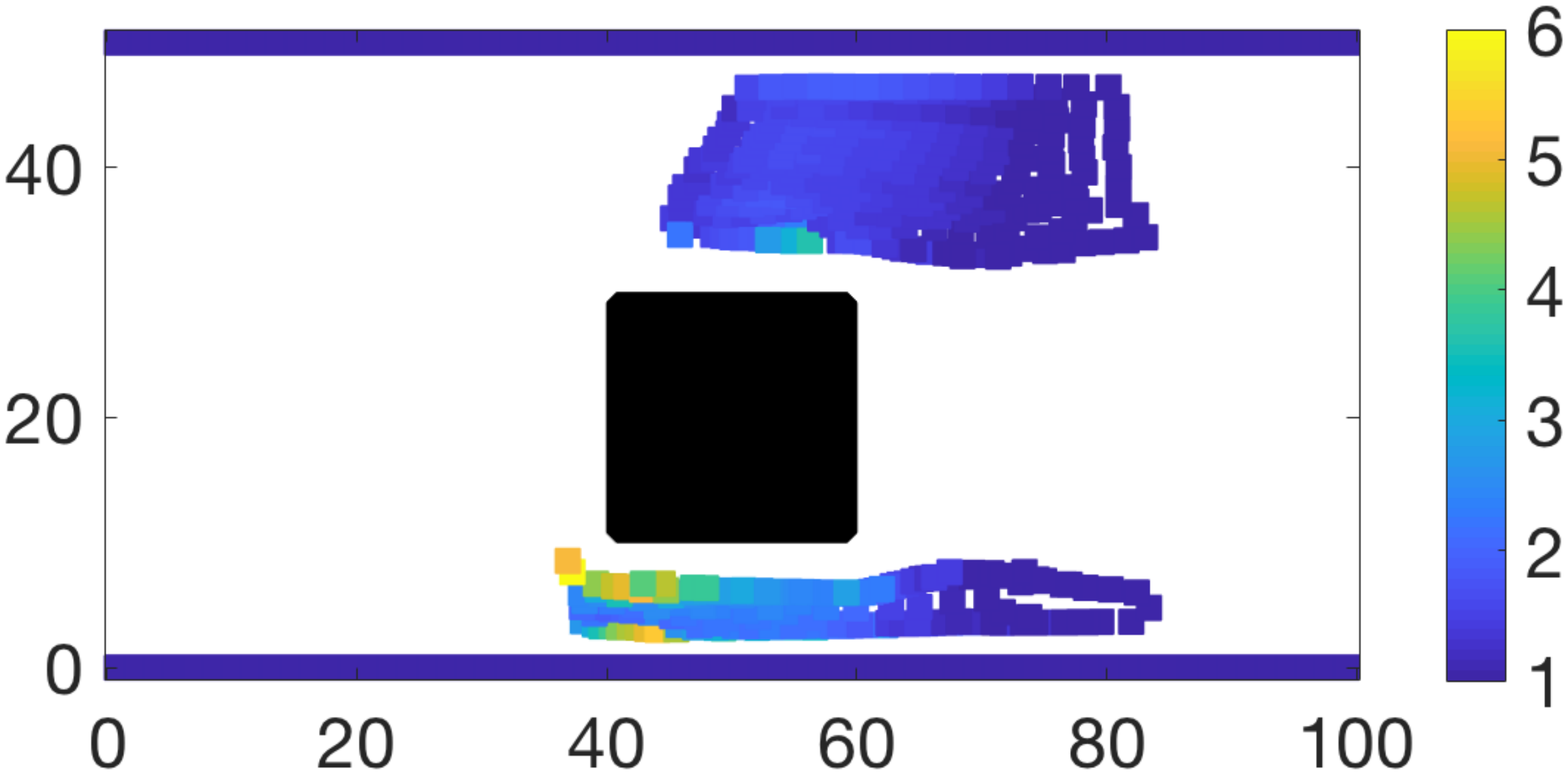}
		\includegraphics[keepaspectratio=true, angle=90, width=0.45\textwidth]{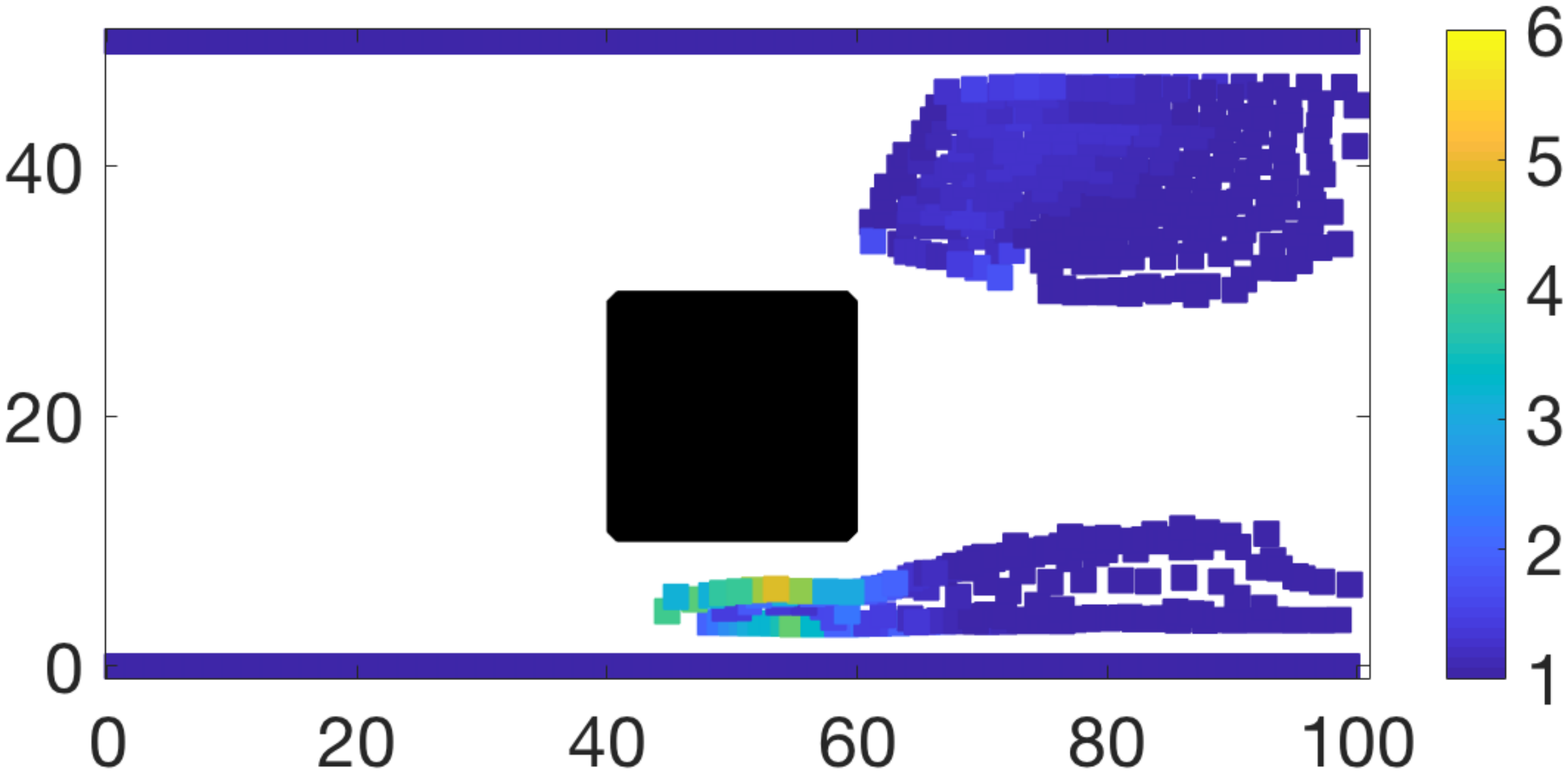}\\		 
		%\vspace{-2cm}
		\caption{ Macroscopic density $\rho$ of pedestrian dynamics  at times $t=10, 20, 30, 40$ with  influence of contact time for uni-direction flow with fixed obstacle. }
		\label{density}
		\centering
	\end{figure}

	Finally, in Figure \ref{density} we show the macroscopic density  of the pedestrians at times $t=10, 20, 30, 40$ for the situation with  influence of contact time for uni-directional  flow with fixed obstacle. 

	\subsection {Test-case 2:  moving obstacle}     
	
	In this subsection we consider the interaction of pedestrians with a  moving obstacle, e.g. a  vehicle in a shared space. We consider  the same computational domain as in test-case 1 with pedestrians, which  are initially located as
	in Figure \ref{phi_v_ne_1_t0} (right bottom). Their destination is the right exit. The moving obstacle of size $4 m \times 8 m$ is located on the right with   the left exit as destination. Typically in a restricted traffic area  the vehicle has
	a low speed limit.  We have chosen  $10 km/h \sim 3m/s$. The maximum speed of the pedestrians is chosen as $2m/s$. 
	
	%\begin{figure}
	%	\centering
	%
	%%	\vspace{-2cm}	 
	%	\caption{Initial situation for pedestrian with moving obstacle at  $t = 0$. Red indicate infected, green indicate susceptibles.Pedestrian dynamics  with  one way (left) and counter flow (right) and without (top) and with obstacle (bottom).}
	%	\label{initial_moving_obstacle}
	%	\centering
	%\end{figure}	     
	
\begin{figure}
		\centering
		\includegraphics[keepaspectratio=true, angle=90, width=0.45\textwidth]{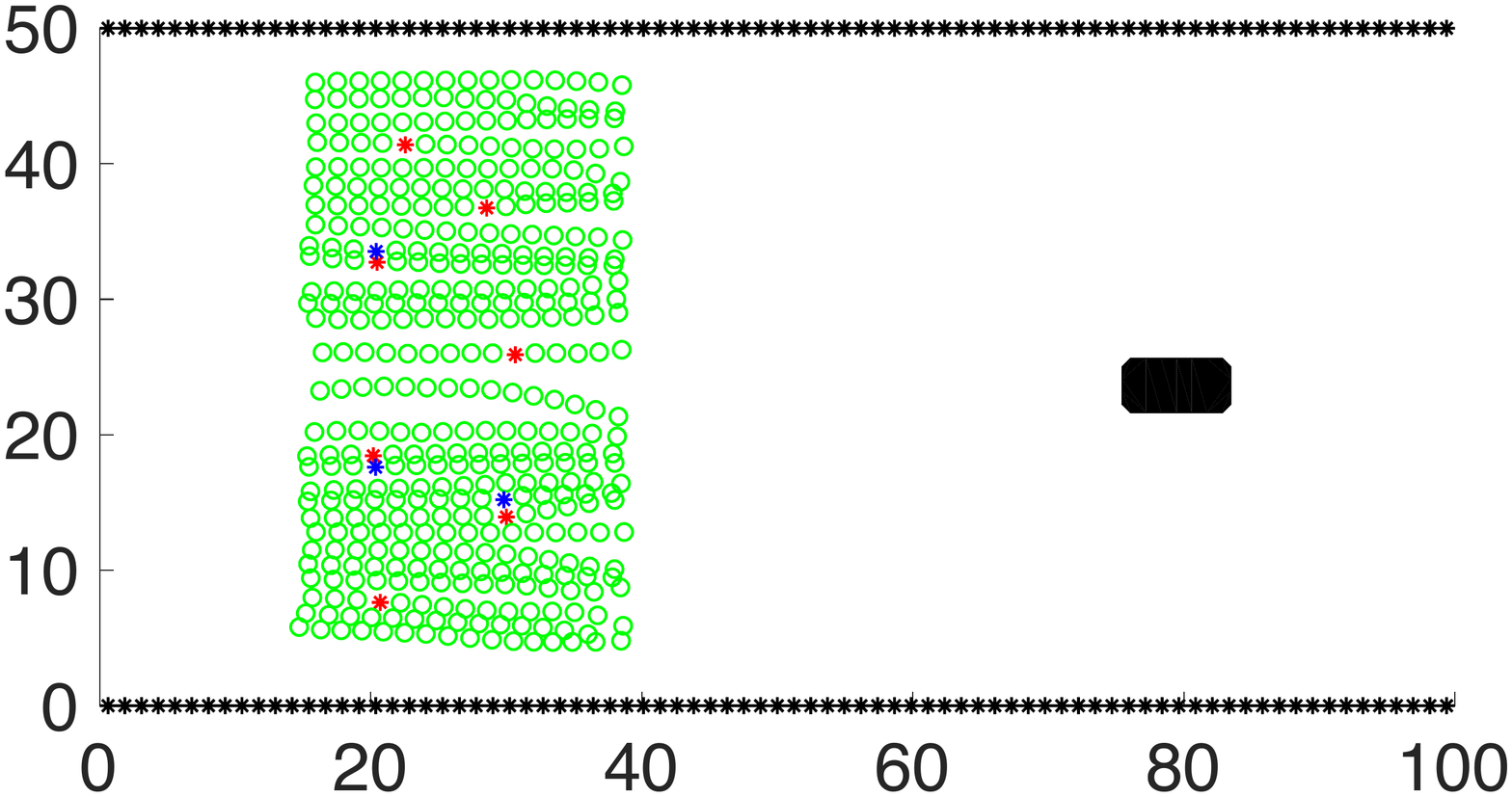}	
		\includegraphics[keepaspectratio=true, angle=90, width=0.45\textwidth]{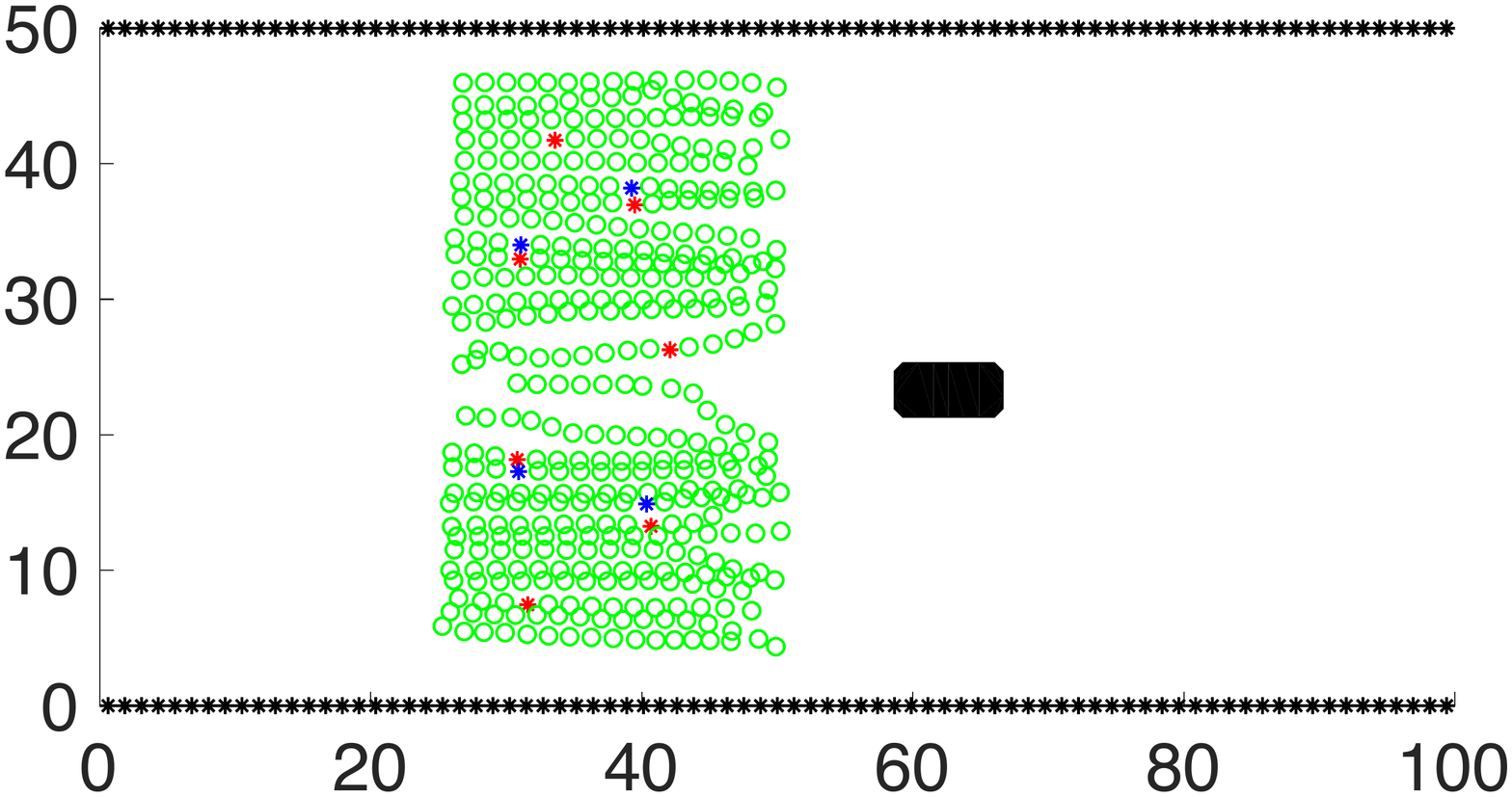}
		\includegraphics[keepaspectratio=true, angle=90, width=0.45\textwidth]{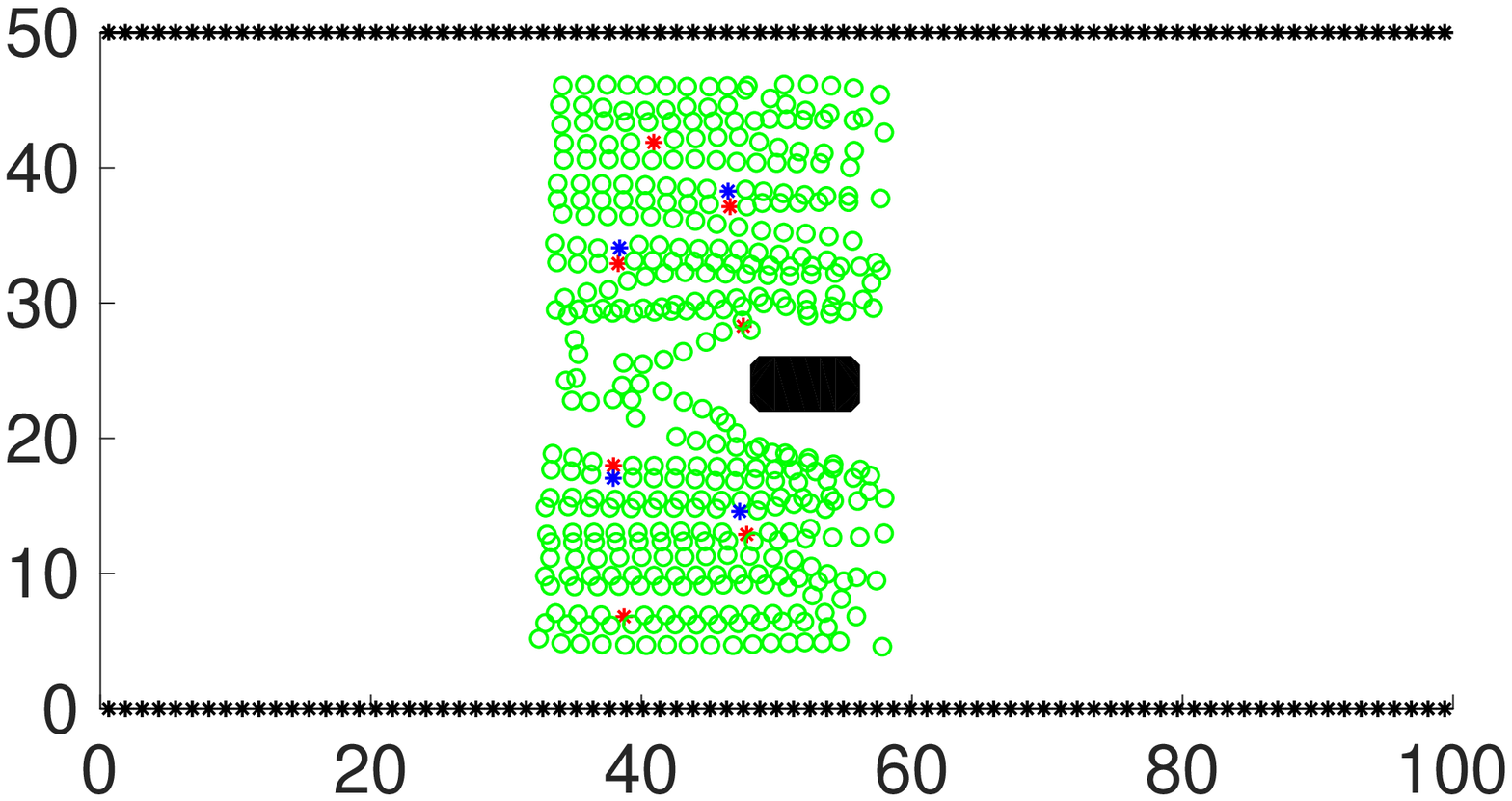}	
		\includegraphics[keepaspectratio=true, angle=90, width=0.45\textwidth]{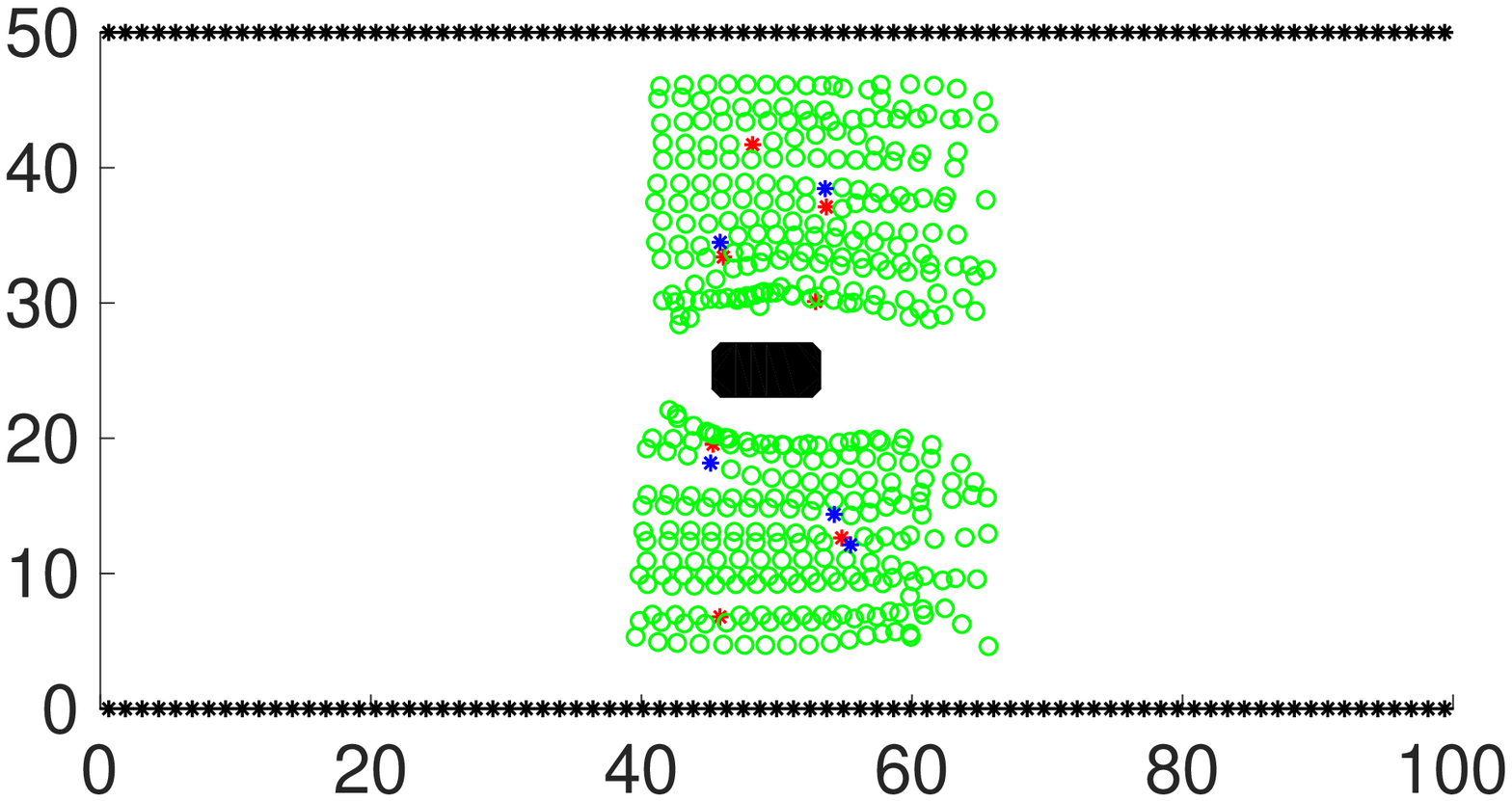}
		\caption{Positions of pedestrian and obstacle at  $t = 8 s, 14 s$ (first row) and $ t = 18 s, 22 s$ (second row). Red indicate infected, green indicate susceptibles, blue indicate probably exposed pedestrians.}
		\label{moving_obstacle}
		\centering
	\end{figure}

	In Figure \ref{moving_obstacle} we have plotted the positions of pedestrians and  obstacle at different times $ t = 8 s, 14 s, 18 s$ and $22 s$. We observe an interaction of pedestrians and obstacle between times $ t = 14 s $ and $t=22s$.  
	Moreover, one observes a slight increase of the number of exposed people which is less pronounced than in the case of the big 
	fixed obstacle in test-case 1.

	% In Figure \ref{CM_obstacle_4} we have plotted the center of mass of the obstacle up to $ t = 39 s$. One observes a  small change of direction of the obstacle, when it interacts with the pedestrians.  
	
	%
	%\begin{figure}
	%	\centering
	%	\includegraphics[keepaspectratio=true, angle=90, width=0.5\textwidth]{figures/fig_with_moving_obstacle/pedestrian_with_moving_obstacle_CM_t_40}	
	%	\caption{Positions of pedestrian and obstacle at time $t=40s$ and  center of mass of the obstacle up to  time $t = 40 s$. Red indicate infected, green indicate susceptibles.}
	%	\label{CM_obstacle_4}
	%	\centering
	%\end{figure}	
	%
	
\begin{figure}
		\centering
		\includegraphics[keepaspectratio=true, angle=90, width=0.6\textwidth]{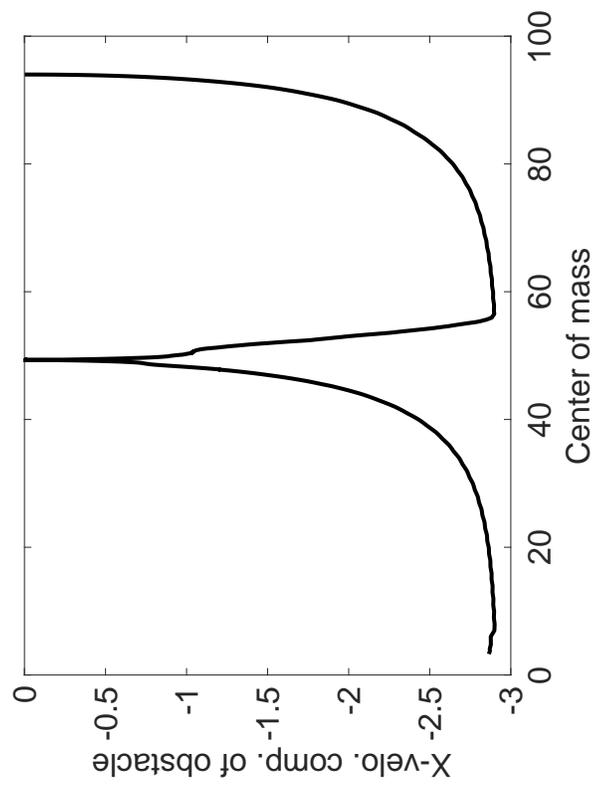}	
		\caption{The velocity of the obstacle along the center of mass. }
		\label{velocity_obstacle}
		\centering
	\end{figure}

	Moreover, in Figure \ref{velocity_obstacle} we have plotted the $x$-velocity  component of the obstacle along its center of mass.  One observes that  the obstacle (coming from the right) accelerates and maintains almost its maximum speed. When it encounters the pedestrian crowd, it reduces  its speed. 
	Finally, it accelerates again, when there are no pedestrians anymore in the surroundings. 
	
	\section{Concluding Remarks}
	We have presented a multi-group macroscopic pedestrian flow model combining a dynamic model for  pedestrians flows and a SEIS based kinetic disease spread model. 
	A meshfree particle method to solve the governing equations is presented and used for the  computation of several numerical examples analyzing  different situations  and parameters.
	The dependence of the solutions and, in particular, the dependence of the number of exposed pedestrians on geometry and parameters is investigated and discussed and shows qualitatively consistent results.
	Findings indicate, that in particular, in bi-directional flow it is important to take into account the contact time for a realistic description of the flow.  This is a realistic qualitative behavior, which sheds a new light for the design of emergency exits in the presence of a pandemic.

	\section*{Acknowledgment}
	
	This work is supported by the German research foundation, DFG grant KL 1105/30-1 and  by the
	DAAD PhD program MIC.

\end{document}